\documentclass[12pt]{amsart}
\usepackage{amsmath}
\usepackage{amssymb}
\newtheorem{theorem}{Theorem}[section]
\newtheorem{lemma}[theorem]{Lemma}
\newtheorem{prop}[theorem]{Proposition}
\newtheorem{cor}[theorem]{Corollary}
\theoremstyle{definition}

\newtheorem{example}[theorem]{Example}

\theoremstyle{remark}
\newtheorem{remark}[theorem]{Remark}
\numberwithin{equation}{section}

\begin{document}

\newcommand\B{\mathbb{B}}
\newcommand\C{\mathbb{C}}
\newcommand\D{\mathbb{D}}
\newcommand\Z{\mathbb{Z}}
\newcommand\Q{\mathbb{Q}}
\newcommand\R{\mathbb{R}}
\newcommand\N{\mathbb{N}}
\newcommand\T{\mathbb{T}}
\newcommand\cD{\mathcal{D}}
\newcommand\cF{\mathcal{F}}
\newcommand\cG{\mathcal{G}}
\newcommand\cH{\mathcal{H}}
\newcommand\cI{\mathcal{I}}
\newcommand\cL{\mathcal{L}}
\newcommand\cM{\mathcal{M}}
\newcommand\cO{\mathcal{O}}
\newcommand\cP{\mathcal{P}}
\newcommand\cS{\mathcal{S}}
\newcommand\cW{\mathcal{W}}
\newcommand\RE{\mathrm{Re}\;}
\newcommand\HD{\cH\cD}
\newcommand\rH{\mathbb{H}_r}
\newcommand\Hi{L^2(0,\infty)}
\newcommand\Ha{H^2(\rH)}
\newcommand\inpr[2]{\langle{#1,#2}\rangle}
\newcommand\Ft[1]{\cF{[#1]}}
\newcommand\Lt[1]{\cL{[#1]}}
\newcommand\Mt[1]{\cM{[#1]}}
\newcommand\Pin[1]{\cP{[#1]}}

\title{A perturbation problem for the shift semigroup}

\author{Masaki Izumi}
\address{Department of Mathematics\\ Graduate School of Science\\
Kyoto University\\ Sakyo-ku, Kyoto 606-8502\\ Japan}
\email{izumi@math.kyoto-u.ac.jp}

\thanks{Work supported by JSPS}

\begin{abstract}
Motivated by B. Tsirelson's construction of 
$E_0$-semigroups of type III, we investigate a $C_0$-semigroup 
acting on the space of square integrable functions of the half line 
whose difference from the shift semigroup is a Hilbert-Schmidt operator. 
We give a description of such semigroups in terms of analytic functions 
on the right-half plane and construct several examples. 
\end{abstract}
\maketitle

\section{Introduction} 
An $E_0$-semigroup is a semigroup of unital endomorphisms of $\B(H)$, 
the set of bounded operators on an infinite dimensional separable 
Hilbert space $H$, with appropriate continuity. 
Despite the fact that $\B(H)$ is the simplest infinite dimensional factor, 
the classification of $E_0$-semigroups is far from a satisfactory stage mainly caused by the 
presence of so called type II and type III examples 
\cite{P1},\cite{P2},\cite{T1}. 
Indeed, these classes are related to various areas of analysis, 
such as probability theory and harmonic analysis, which makes the subject 
more interesting and worth investigating than one would expect at first sight. 
The reader is referred to Arveson's monograph \cite{Arv} and 
contributions in \cite{Pr}, in particular those of Arveson, Powers, 
and Tsirelson among others, for history and basic results for the subject. 

In this paper, we investigate a purely operator theoretical problem 
motivated by Tsirelson's construction of uncountably many mutually 
non-cocycle conjugate $E_0$-semigroups of type III \cite{T2}. 
Namely, let $\{S_t\}_{t\geq 0}$ be the shift semigroup of 
$\Hi$, that is
$$(S_tf)(x)=\left\{
\begin{array}{ll}0& (x<t)\\
f(x-t) &(t\leq x)
\end{array}
\right. .$$
We investigate the structure of a $C_0$-semigroup of bounded operators 
$\{T_t\}_{t\geq 0}$ acting on $\Hi$ satisfying the following two conditions: 
\begin{itemize}
\item [(C1)] $T_t^*S_t=I$ for all $t\geq 0$. 
\item [(C2)] $T_t-S_t$ is a Hilbert-Schmidt class operator for all 
$t\geq 0$.   
\end{itemize}
If $T_t$ as above preserves the real functions $\Hi_{\R}$ in $\Hi$,  
then a standard argument for quasi-equivalence of representations of the  
CCR (canonical commutation relation) algebra tells the following: 
there exists an $E_0$-semigroup $\alpha_t$ acting on the symmetric Fock space 
over $\Hi$ determined by 
$$\alpha_t(W(f+ig))=W(S_tf+iT_tg),\quad f,g\in \Hi_{\R},$$ 
where $W(f+ig)$ is the Weyl operator for $f+ig$. 
Although Tsirelson's construction is based on a homogeneous continuous 
product of measure classes arising from an off white noise \cite{T1}, 
his examples actually come from the above construction as we will see 
in Section 6. 

A careful examination of the generator of a $C_0$-semigroup 
$\{T_t\}_{t\geq 0}$ satisfying (C1) shows that $\{T_t\}_{t\geq 0}$ 
is completely characterized  by a holomorphic function $M(z)$ on the 
right half plane $\rH$ such that $M(z)/(1+z)$ belongs to the Hardy space 
$\Ha$. 
We call $M(z)$ the half-density function for $\{T_t\}_{t\geq 0}$. 
When $\{T_t\}_{t\geq 0}$ comes from an off white noise, the absolute value of 
the boundary value function $M(iy)$ on the imaginary axis is nothing but 
the square root of the spectral density function discussed by Tsirelson 
\cite{T1},\cite{T2}. 
We completely characterize $C_0$-semigroups satisfying (C1) and (C2) in terms 
of the function $M(z)$ (Theorem 4.2). 
It turns out that this class of functions $|M(i\lambda)|^2$ is slightly larger than the class 
of spectral density functions discussed in \cite{T2}. 

In Section 5, we construct a class of $C_0$-semigroups satisfying (C1) and 
(C2) parameterized by a function $\varphi$ in the space 
$L^1_{\mathrm{loc}}[0,\infty)\cap L^2((0,\infty),(1\wedge x)dx)$.  
More precisely, we consider the case with $M(z)=1-\Lt \varphi(z)$ where 
$\Lt \varphi$ is the Laplace transform of $\varphi$. 
While a perturbation argument shows that the $E_0$-semigroups arising from 
this class of functions are cocycle conjugate to those coming from a subclass of off white noises, 
one cannot distinguish them from the CCR flow of index 1 by using 
Tsirelson's invariants discussed in \cite{T1} because the spectral density 
functions $|M(iy)|^2$ converge to 1 at infinity in this case. 
Yet, in the forthcoming paper \cite{IS}, we show that this class of 
$E_0$-semigroups contains uncountably many mutually non cocycle conjugate 
type III $E_0$-semigroups. 
The invariant we adopt for differentiating these examples is 
the Murray-von Neumann type of the ``local observable algebras" for an open 
subset $U$ of the interval $[0,1]$, which may be an AFD type III factor if 
$\varphi$ does not belong to $L^2(0,\infty)$ and $U$ has a sufficiently 
complicated shape. 
It should be noted that the space $L^2((0,\infty),(1\wedge x)dx)$ also 
plays a crucial role in 
Powers' construction of uncountably many mutually non cocycle conjugate 
$E_0$-semigroups of type II$_0$ \cite{P2} though we don't known if 
there exists a direct relationship between our argument with Powers' 
in the present stage. 

The author would like to thank W. Arveson and R. Srinivasan for 
stimulating discussions.

\section{Preliminaries} 
Throughout this note, the symbol $||\cdot||$ means either the $L^2$ norm or the 
operator norm, depending on the context. 
The symbol $||\cdot ||_p$ denotes the $L^p$ norm. 
The Hilbert-Schmidt norm is denoted by $||\cdot||_{H.S.}$. 
Every function space we discuss is complex valued. 

For $f,g\in \Hi$, we denote by $\inpr{f}{g}$ the usual inner product 
$$\inpr{f}{g}=\int_0^\infty f(x)\overline{g(x)}dx.$$
Dual pairing in various contexts will be denoted by $(f,g)$. 
For example, for $f,g\in \Hi$ 
$$(f,g)=\int_0^\infty f(x)g(x)dx.$$

For an operator $A$ of $\Hi$, we denote by $D(A)$ the domain of $A$ 
and  by $\cG(A)$ the graph of $A$, that is  
$$\cG(A)=\{f\oplus Af;\; f\in D(A)\}.$$
We set 
$$\cG'(A)=\{-Af\oplus f;\; f\in D(A)\}.$$
It is well-known that when $A$ is densely defined closable 
operator, we have $\cG(A)^\perp=\cG'(A^*)$. 
We sometimes use the notation 
$$\inpr{f}{g}_A=\inpr{f}{g}+\inpr{Af}{Ag}$$
$$(f,g)_A=(f,g)+(Af,Ag).$$

For $f,g\in \Hi$, we define $f\otimes g\in \B(\Hi)$ by 
$$(f\otimes g)h=(h,g)f,\quad h\in \Hi.$$

We denote by $\rH$ the right-half plane $\{z\in \C;\; \RE z>0\}$. 
For $z\in \rH$, we set $e_z(x)=e^{-zx}$. 
We denote by $L^1_{\mathrm{loc}}[0,\infty)$ the set of measurable 
functions on $[0,\infty)$ that are integrable on every compact subset of 
$[0,\infty)$. 
For $f\in L^1_{\mathrm{loc}}[0,\infty)$ with $a>1$ such that 
$e_af\in L^1(0,\infty)$, we denote by $\Lt h(z)$ 
the Laplace transformation 
$$\Lt f(z)=\int_0^\infty f(x)e^{-zx}dx,\quad \RE z>a.$$
For $f\in \Hi$, we have $\Lt f(z)=(f,e_z)$ and 
$|\Lt f(z)|\leq ||fe_{z/2}||\cdot||e_{z/2}||$, and so

\begin{equation}
\sqrt{x}|\Lt f(x+iy)|\leq ||fe_{x/2}||\rightarrow 0,\quad 
(x\to+\infty).
\end{equation} 

Let $S=\{S_t\}_{t\geq 0}$ be the shift semigroup of $\Hi$, 
that is
$$(S_tf)(x)=\left\{
\begin{array}{ll}0& (x<t)\\
f(x-t) &(t\leq x)
\end{array}
\right..$$

The main purpose of this paper is to 
investigate the structure of $C_0$-semigroups $T=\{T_t\}_{t\geq 0}$ of 
bounded operators in $\B(\Hi)$ satisfying the following two conditions:
\begin{itemize}
\item [(C1)] $T_t^*S_t=I$ for all $t\geq 0$. 
\item [(C2)] $T_t-S_t$ is a Hilbert-Schmidt class operator for all 
$t\geq 0$.   
\end{itemize}
 
Whenever a function $f\in \Hi$ extends to a continuous function 
on $[0,\infty)$, we abuse the notation and use the same symbol $f$ 
for the extension. 
The generator $A$ of $\{S_t\}_{t\geq 0}$ is the differential operator 
$(Af)(x)=-f'(x)$ with boundary condition $f(0)=0$. 
More precisely, the domain $D(A)$ of $A$ consists of square integrable 
locally absolutely continuous functions $f$ on $[0,\infty)$ such that $f(0)=0$ 
and the derivative $f'$ (well-defined almost everywhere) belongs to $\Hi$. 
The adjoint operator $A^*$, which is the generator of 
$\{S_t^*\}_{t\geq 0}$, is the differential operator $(A^*f)(x)=f'(x)$ 
without any boundary condition. 
More precisely, $D(A^*)$ is the linear span of $D(A)$ and $e^{-x}$.

Let $B$ be the generator of a semigroup $\{T_t\}_{t\geq 0}$ satisfying 
the condition (C1). 
Then for $f\in D(A)$ and $g\in D(B)$ we have 
$\inpr{S_tf}{T_tg}=\inpr{f}{g}$ and  
$$\inpr{Af}{g}+\inpr{f}{Bg}=0,$$
which shows $B\subset -A^*$. 
This means that $B$ is also a differential operator, but what matters now is 
the domain of $B$. 

For $p\in D(A^*)\setminus \{0\}$, we define $A_p\subset -A^*$ whose graph is 
$\cG(-A^*)\cap\C(\overline{p}\oplus -\overline{p'})^\perp$, 
where $\overline{p}$ is the complex conjugate of $p$. 
By definition, $A_p$ is a closed operator. 
As the codimension of $\cG(A_p)$ in $\cG(-A^*)$ is 1, the operator 
$B$ is a restriction of $A_p$ for some $p$ (in fact, we will show 
$B=A_p$ later). 
Note that when $p(x)=e^{-x}$, the operator $A_p$ is nothing but $A$. 

\begin{lemma} The operator $A_p$ is densely defined if and only if 
$p'\not\in D(A)\setminus \{0\}$. 
\end{lemma}

\begin{proof} Note that 
$$D(A_p)^\perp\oplus 0
=(\Hi\oplus 0)\cap \cG(A_p)^\perp=
(\Hi\oplus 0)\cap \Big(\cG'(-A)+
\C(\overline{p}\oplus -\overline{p'})\Big).$$
Thus $g\in D(A_p)^\perp$ if and only if there exist $f\in D(A)$ and 
$\lambda\in \C$ such that 
$$g=-f'+\lambda \overline{p},$$
$$0=f-\lambda \overline{p'}.$$
This shows that if $p'\not\in D(A)$, $g=0$ and so 
the operator $A_p$ is densely defined. 

On the other hand, assume that $p'\in D(A)\setminus {0}$. 
Then $\lambda\overline{(p-p'')}$ belongs to $D(A_p)^\perp$. 
Note that $p-p''\neq 0$ as far as $p'\in D(A)\setminus \{0\}$. 
\end{proof} 

We set 
$$\cO=\{p\in D(A^*);\; p'\notin D(A)\}.$$

\begin{cor} Let $\{T_t\}_{t\geq 0}$ be a $C_0$-semigroup satisfying (C1) 
and let $B$ be the generator of $\{T_t\}_{t\geq 0}$. 
Then there exists $p\in \cO$ such that $B\subset A_p$. 
\end{cor}

We end this section with giving a proof for the claim stated in 
Introduction that the semigroup $\{T_t\}_{t\geq 0}$ satisfying (C1) and (C2) 
gives rise to an $E_0$-semigroup. 
Let $H=e^{L^2(0,\infty)}$ be the symmetric Fock space for $L^2(0,\infty)$ 
and let $\Omega$ be the vacuum vector. 
We denote by $\exp(f)\in H$ the exponential vector for $f\in \Hi$, 
that is 
$$\exp(f)=\bigoplus_{n=0}^\infty \frac{1}{\sqrt{n!}}f^{\otimes^n},$$
where $f^{\otimes^0}=\Omega$. 
The Weyl operator $W(h)$ for $h\in \Hi$ is a unitary operator in $\B(H)$ 
determined by 
$$W(h)\exp(f)=e^{-\frac{\|h\|^2}{2}-\inpr{f}{h}}\exp(f+h).$$
The CCR algebra for $\Hi$ is the algebra generated by $\{W(h)\}_{h\in \Hi}$ 
whose defining (vacuum) representation is irreducible. 

\begin{theorem} If $\{T_t\}_{t\geq 0}$ satisfies (C1) and (C2) and $T_t$ 
preserves the real functions $\Hi_\R$, then there exists a unique 
$E_0$-semigroup $\{\alpha\}_{t\geq 0}$ acting on $\B(H)$ satisfying.
$$\alpha_t(W(f+ig))=W(S_tf+iT_tg),\quad f,g\in \Hi_{\R}.$$ 
\end{theorem}

\begin{proof} The condition (C1) shows that there exists a representation 
$\beta_t$ of the CCR algebra for $\Hi$ satisfying 
$\beta_t(W(f+ig))=W(S_tf+iT_tg)$ for 
$f,g\in \Hi_\R$. 
To show the claim, it suffices to prove that $\beta_t$ is quasi-equivalent 
to the vacuum representation. 
We set $M$ to be the von Neumann algebra generated by the image of $\beta_t$. 
Since $M$ is a factor thanks to \cite[Theorem 1]{A1}, it suffices to show that 
the representation $\beta_t$ restricted to $\overline{M\Omega}$ is 
quasi-equivalent to the vacuum representation. 
Thus our problem is now reduced to the quasi-equivalence of the vacuum 
representation and the GNS representation of the quasi-free state 
$\omega'(x)=\inpr {\beta_t(x)\Omega}{\Omega}$, which can be shown 
by using (C2) from a well-known criterion (see \cite{A2},\cite{vD} 
for example). 
\end{proof}

\section{Resolvent}
In this section, we explicitly compute the resolvent of the generator of 
a semigroup $\{T_t\}_{t\geq 0}$ satisfying (C1). 
It turns our that such a semigroup is completely characterized by 
a holomorphic function on the right-half plane. 

For $p\in D(A^*)$, we set 
\begin{eqnarray*}\Mt p(z)&:=&(p,e_z)_{A^*}=\Lt p(z)-z\Lt{p'}(z)=(1-z^2)\Lt p(z)+p(0)z\\
&=&(\frac{1}{z}-z)\Lt {p'}(z)+\frac{p(0)}{z}.
\end{eqnarray*}
Note that unless $p=0$, the function $\Mt p(z)$ is not constantly zero. 

For a point $z\in \rH$ with $\Mt p(z)\neq 0$, we define $\xi_{p,z}(x)$ by
\begin{eqnarray*}\xi_{p,z}(x)&=&-\frac{\int_0^\infty e^{-zt}\Big(p(x+t)
+z(p'(x)-p'(x+t))\Big)dt}{\Mt p(z)}\\
&=&-\frac{(1-z^2)\int_0^\infty e^{-zt}p(x+t)dt+zp(x)+p'(x)}{\Mt p(z)}.
\end{eqnarray*}
Note that $\xi_{p,z}$ is given by the $\Hi$-valued 
integral 
$$\xi_{p,z}=-\frac{(1-z^2)\int_0^\infty e^{-zt}S_t^*pdt+zp+p'}{\Mt p(z)}$$
and so $\xi_{p,z}\in \Hi$. 

\begin{lemma} Let $p\in \cO$ and $\RE z>0$. 
Then $z$ is in the resolvent set of $A_p$ if and only if 
$\Mt p(z)\neq 0$. 
If $z$ satisfies this condition, the resolvent $(zI-A_p)^{-1}$ 
is given by 
$$(zI-A_p)^{-1}=(zI-A)^{-1}+e_z\otimes\xi_{p,z}.$$ 
\end{lemma}

\begin{proof} If $\Mt p(z)=0$, we have $e_z\in D(A_p)$ and $A_pe_z=ze_z$. 
Thus $\Mt p(z)\neq 0$ is a necessary condition for $z$ to 
belong to the resolvent set of $A_p$.  

Assume $\Mt p(z)\neq 0$ now. 
Let 
$$B_z=(zI-A)^{-1}+e_z\otimes\xi_{p,z}.$$
Note that $z$ belongs to the 
resolvent set of $A$ for $\RE z>0$ and the resolvent is given by 
$$(zI-A)^{-1}f(x)=\int_0^xf(t)e^{-z(x-t)}dt.$$
First we show $B_z(zI-A_p)g=g$ for $g\in D(A_p)$. 
As 
$$(zI-A)^{-1}(zI-A_p)g(x)=\int_0^x(zg(t)+g'(t))e^{-z(x-t)}dt
=g(x)-g(0)e^{-zx},$$
it suffices to show $(zI-A_pg,\xi_{p,z})=g(0)$. 
Indeed, 
\begin{eqnarray*}((zI-A_p)g,\xi_{p,z})&=&
-\frac{\int_0^\infty e^{-zt}((zI-A_p)g,S_t^*(p-zp'))dt+(zg+g',p')}
{\Mt p(z)}\\
&=&-\frac{\int_0^\infty e^{-zt}(S_t(zI-A_p)g,p-zp')dt+(zg+g',p')}
{\Mt p(z)}\\
&=&-\frac{((zI-A)^{-1}(zI-A_p)g,p-zp')+(zg+g',p')}
{\Mt p(z)}\\
&=&-\frac{(g-g(0)e_z,p-zp')+(zg+g',p')}
{\Mt p(z)}\\
&=&-\frac{(g,p)+(g',p')-g(0)\Big((e_z,p)-z(e_z,p')\Big)}
{\Mt p(z)}\\
&=&g(0). 
\end{eqnarray*}

Next we show that the range of $B_z$ is contained in $D(A_p)$. 
Note that for $f\in \Hi$ the element $B_zf=(zI-A)^{-1}f+(f,\xi_{p,z})e_z$ 
belongs to $D(A^*)$, and so 
\begin{eqnarray*}\lefteqn{(B_zf,p)+(A^*B_zf,A^*p)}\\
&=&((zI-A)^{-1}f,p)+(A^*(zI-A)^{-1}f,A^*p)
+\Big((e_z,p)+(A^*e_z,A^*p)\Big)(f,\xi_{p,z})\\
&=&((zI-A)^{-1}f,p)-(A(zI-A)^{-1}f,A^*p)+\Mt p(z)(f,\xi_{p,z})\\
&=&((zI-A)^{-1}f,p)+(f-z(zI-A)^{-1}f,A^*p)+\Mt p(z)(f,\xi_{p,z})\\
&=&((zI-A)^{-1}f,p)+(f,A^*p)-z(A(zI-A)^{-1}f,p)+\Mt p(z)(f,\xi_{p,z})\\
&=&(1-z^2)((zI-A)^{-1}f,p)+(f,p')+z(f,p)+\Mt p(z)(f,\xi_{p,z})\\
&=&(1-z^2)((zI-A)^{-1}f,p)-(1-z^2)\int_0^\infty e^{-zt}(f,S_t^*p)\\
&=&(1-z^2)((zI-A)^{-1}f,p)-(1-z^2)\int_0^\infty e^{-zt}(S_tf,p)\\
&=&0. 
\end{eqnarray*} 
Thus, we get $B_zf\in D(A_p)$ and  
\begin{eqnarray*}(zI-A_p)B_zf&=&
(zI+A^*)\Big((zI-A)^{-1}f+(f,\xi_{p,z})e_z\Big)
\\&=&(zI+A^*)(zI-A)^{-1}f=f.
\end{eqnarray*}
Therefore, $B_z=(zI-A_p)^{-1}$.
\end{proof}

Note that when $A_1$ and $A_2$ are densely defined closed operators 
on a Banach space such that $A_1$ is a proper extension of $A_2$, 
the intersection of the resolvent sets of $A_1$ and $A_2$ is empty.

\begin{theorem} Let $\{T_t\}_{t\geq 0}$ be a $C_0$-semigroup 
satisfying (C1). 
Then there exists $p\in \cO$ such that 
$A_p$ is the generator of $\{T_t\}_{t\geq 0}$. 
Moreover, there exists a positive number $a$ such that 
for $\forall f\in D(A_p)$, $\forall t>0$ and $\forall b>a$, 
the function $T_tf$ is given by 
$$T_tf=S_tf+\lim_{r\to+\infty}\frac{1}{2\pi i}\int_{b-ir}^{b+ir}
e^{zt}(f,\xi_{p,z})e_zdz.$$
\end{theorem}

\begin{proof} Let $B$ be a generator of $\{T_t\}_{t\geq 0}$. 
Thanks to Corollary 2.2, there exists $p\in D(A^*)$ with 
$p'\not\in D(A)$ such that $B\subset A_p$. 
It is well-known that there exist positive numbers $a$ and $M$ such that 
for all $t\geq 0$ the inequality $||T_t||\leq Me^{at}$  holds \cite[p.232]{Y}.
Thus $z$ belongs to the resolvent set of $B$ for $\RE z>a$ and 
$$(zI-B)^{-1}f=\int_0^\infty e^{-zt}T_tf,\quad \forall f\in \Hi.$$
On the other hand, since $\Mt p(z)$ is a non-constant analytic function on 
$\rH$, there exists $z_0$ such that $\RE z_0>a$ and $\Mt p(z_0)\neq 0$, 
and so Lemma 3.2 implies that $z_0$ is in the resolvent set of $A_p$. 
If $B$ were proper restriction of $A_p$, the point $z_0$ would not belong to 
the resolvent set of $B$, which is a contradiction. 
Thus $B=A_p$. 
The integral formula above follows from the usual inverse Laplace 
transformation with the fact that the map $(0,\infty)\ni t\mapsto T_tf\in 
\Hi$ is differentiable for $f\in D(A_p)$.    
\end{proof}

The above theorem shows that every information of the semigroup 
$\{T_t\}_{t\geq 0}$ is encoded in $\xi_{p,z}$, which is determined 
by $\Mt p$ via the following lemma: 

\begin{lemma} Let $p\in \cO$ and let $z$ satisfy $\RE z>0$ and 
$\Mt p(z)\neq 0$. 
Then 
$$\Lt{\xi_{p,z}}(w)=\frac{\Mt p(z)-\Mt p(w)}{\Mt p(z)(z-w)}.$$
\end{lemma}
\begin{proof} Assume $w\neq z$. 
Then
\begin{eqnarray*}\lefteqn{\Mt p(z)\Lt{\xi_{p,z}}(w)}\\
&=&
-\Lt{p'}(w)-z\Lt p(w)-(1-z^2)\int_0^\infty dt
\int_0^\infty dx e^{-wx}e^{-zt}p(x+t)\\
&=&p(0)-(z+w)\Lt p(w)-(1-z^2)\int_0^\infty dt
\int_t^\infty dx e^{-w(x-t)}e^{-zt}p(x)\\
&=&p(0)-(z+w)\Lt p(w)-(1-z^2)\int_0^\infty p(x)e^{-wx}dx
\int_0^x dte^{-(z-w)t}\\
&=&p(0)-(z+w)\Lt p(w)+(1-z^2)\int_0^\infty p(x)
\frac{e^{-zx}-e^{-wx}}{z-w}dx\\
&=&p(0)-(z+w)\Lt p(w)+(1-z^2)\frac{\Lt p(z)-\Lt p(w)}{z-w}\\
&=&\frac{\Mt p(z)-\Mt p(w)}{z-w},\\
\end{eqnarray*}
which shows the statement. 
\end{proof}

When $A_p$ generates a $C_0$-semigroup $\{T_t\}_{t>0}$, 
we call $p$ the {\it orthogonal function} of $\{T_t\}_{t>0}$ and call $\Mt p$ 
the {\it half density function} of $\{T_t\}_{t>0}$. 
The half density function plays the central role throughout 
this note and the theory of the Hardy classes is best suit for its 
description. 

We denote by $H^p(\rH)$ the Hardy space of the right half plane $\rH$, 
that is, the set of holomorphic functions $F$ on $\rH$ 
such that $\sup_{x>0}||F(x+i\cdot)||_p$ is finite. 
For $F\in H^p(\rH)$ and $iy$ in the imaginary axis, 
we denote by $F(iy)$ the non-tangential limit of $F$ at $iy$, which makes 
sense for almost every $iy$. 
Note that when $f\in \Hi$, its Laplace transform $\Lt f$ belongs to 
$\Ha$. 
Thanks to the Paley-Wiener theorem, every element of $\Ha$ is of this form. 

Let $\HD$ be the set of holomorphic functions $M(z)$ 
on $\rH$ such that $M(z)/(1+z)$ belongs to $\Ha$ and 
$M$ does not belong to $\Ha$. 

\begin{lemma} Let the notation be as above. 
Then,
\begin{itemize}
\item [$(1)$] The map $p\mapsto \Mt p$ gives a bijection between $\cO$ and 
$\HD$. 
\item [$(2)$] The map $q\mapsto (1+z)\Lt q(z)$ 
gives a bijection between $\Hi\setminus D(A)$ and $\HD$. 
\item [$(3)$] The bijection between $\cO$ and $\Hi\setminus D(A)$ 
induced by $(1)$ and $(2)$ is given by 
$$q(x)=p(x)-p'(x)-p(0)e^{-x},$$
$$p(x)=2\Lt q(1)\cosh x-\int_0^x e^{x-t}q(t)dt
=\Lt q(1)e^{-x}+\int_0^\infty e^{-t}S_t^*q(x)dt.$$
\end{itemize}
\end{lemma}

\begin{proof} (2)
Since $D(A)$ is exactly the set of $q\in \Hi$ with $z\Lt q(z)\in \Ha$, 
this is a direct consequence of the the Paley-Winer theorem. 

(1) Let $p\in \cO$. 
Using $\Lt {p'}(x)=z\Lt p(z)-p(0)$, we get 
$$\frac{\Mt p(z)}{1+z}=\Lt p(z)-\Lt{p'}(z)-\frac{p(0)}{1+z}
=\Lt{p-p'-p(0)e_1}(z)\in \Ha.$$ 
Assume that $\Mt p(z)$ belongs to $\Ha$. 
Then since we have
$$\frac{z\Mt p(z)}{1+z}=(1-z)\Lt{p'}(z)+\frac{p(0)}{1+z},$$
we would get 
$$z\Lt{p'}(z)=\Lt{p'}(z)-\frac{z\Mt p(z)}{1+z}+\frac{p(0)}{1+z}\in \Ha.$$
This implies $p'\in D(A)$, which is contradiction. 
Thus $\Mt p\in \HD$. 

Injectivity of the map $p\mapsto \Mt p$ is obvious and we show 
surjectivity now. 
For $M\in \HD$, take $q\in \Hi\setminus D(A)$ satisfying 
$M(z)=(1+z)\Lt q(z)$. 
We define $p(x)$ by the second formula in (3). 
Then $p$ is an absolutely continuous function in $\Hi$ satisfying 
the first formula and so $p\in D(A^*)$. 
Since $p'-q=p-p(0)e_1\in D(A)$ and $g$ does not belong to $D(A)$, 
we conclude $p'\notin D(A)$ and so $p\in \cO$. 
$M=\Mt p$ can be shown by computing the Laplace transform. 

(3) has already been shown in the above argument.  
\end{proof}

From now on, we are mainly working on $\HD$ and $\Hi\setminus D(A)$ 
instead of $\cO$. 
When $p\in \cO$ and $M=\Mt p$, we abuse the notation and use the symbols 
$A_M$ and $\xi_{M,z}$ for $A_p$ and $\xi_{p,z}$ respectively. 
Lemma 3.3 means  
\begin{equation}
\Lt{\xi_{M,z}}(w)=\frac{M(z)-M(w)}{M(z)(z-w)}.
\end{equation}
When $(1+z)\Lt q(z)=M(z)$, we can show the following 
by computing the Laplace transform:
\begin{equation}
\xi_{M,z}(y)=\frac{q(y)-(1+z)\int_0^\infty q(y+s)e^{-sz}ds}{M(z)}. 
\end{equation}

We denote by $\HD_b$ the set of $M\in \HD$ such that $A_M$ generates a 
$C_0$-semigroup $\{T_t\}_{t>0}$. 
Such $\{T_t\}_{t>0}$ always satisfies (C1). 
We denote by $\HD_2$ the set of $M\in \HD_b$ such that $\{T_t\}_{t>0}$ 
satisfies (C1) and (C2). 

When $x+iy\in \rH$ and $f$ is a measurable function on $\R$ such that 
$f(\lambda)/(1+\lambda^2)\in L^1(\R)$, we denote by $\Pin f(x+iy)$ 
the Poisson integral 
$$\Pin f(x+iy)=\frac{1}{\pi}\int_{-\infty}^\infty\frac{xf(\lambda)}
{x^2+(y-\lambda)^2}d\lambda.$$

\begin{lemma} Let $M\in \HD$ and assume $z=x+iy\in \rH$ is not a zero of 
$M$. 
Then
$$||\xi_{M,z}||^2=\frac{1}{2\RE z}
\big(\frac{\Pin{|M(i\cdot)|^2}(z)}{|M(z)|^2}-1\big).$$
\end{lemma}

\begin{proof} 
The Plancherel theorem implies 
\begin{eqnarray*}||\xi_{M,z}||^2
&=&\frac{1}{2\pi}
\int_{-\infty }^{+\infty}
|\Lt {\xi_{M,z}}(i\lambda)|^2d\lambda\\
&=&\frac{1}{2\pi}
\int_{-\infty}^{+\infty}
\frac{|M(z)-M(i\lambda)|^2}{|M(z)(z-i\lambda)|^2}d\lambda\\
&=&\frac{1}{2\pi}
\int_{-\infty}^{+\infty}
\frac{|M(z)|^2-2\RE\big(M(i\lambda)\overline{M(z)}\big)
+|M(i\lambda)|^2}{|M(z)|^2|z-i\lambda|^2}d\lambda\\
&=&\frac{1}{2\RE z|M(z)|^2}\big[|M(z)|^2-
2\RE\big(\Pin{M(i\cdot)}(z)\overline{M(z)}\big)
+\Pin{|M(i\cdot)|^2}(z)\big].
\end{eqnarray*}
Thus it suffices to show $\Pin{M(i\cdot)}(z)=M(z)$. 
Let $re^{i\theta}=(z-1)/(z+1)$ and set 
$f(re^{i\theta})=M(z)$. 
Since $M(z)/(1+z)\in \Ha$, the function $f$ belongs to the Hardy space 
$H^2(\D)$ of the unit disc $\D$ \cite[Chapter 8]{H} and 
$$\Pin{M(i\cdot)}(z)=\frac{1}{2\pi}\int_{0}^{2\pi}
\frac{(1-r^2)f(e^{it})}{1-2r\cos(\theta-t)+r^2}dt.$$
Therefore $\Pin{M(i\cdot)}(z)=M(z)$ follows from the 
usual Poisson integral formula. 
\end{proof}

Note that $\Pin{|M(i\cdot)|^2}(z)\geq |M(z)|^2$ always holds.

\begin{cor} Let $M\in \HD$. 
If $A_M$ generates a $C_0$-semigroup, then 
there exist positive constants $a$ and $C$ such that for all 
$\RE z\geq a$, 
$$\log\Pin{|M(i\cdot)|^2}(z)-\log|M(z)|^2\leq 
\frac{\Pin{|M(i\cdot)|^2}(z)}{|M(z)|^2}-1\leq C. $$
Moreover, 
$$\sqrt{x}|M(x)|\rightarrow +\infty\quad (x\to +\infty).$$
\end{cor}

\begin{proof}
Lemma 3.1 shows that if $M\in \HD$ and $M(z)\neq 0$, 
\begin{eqnarray*}
\frac{||\xi_{M,z}||}{\sqrt{2\RE z}}&=&||e_z\otimes \xi_{M,z}||
=||(zI-A_h)^{-1}-(zI-A)^{-1}||\\
&\leq& ||(zI-A_h)^{-1}||+||(zI-A)^{-1}||.
\end{eqnarray*}
When $A_M$ generates a $C_0$-semigroup, there exist constant $a>0$ 
and $C_1$ such that $||(zI-A_M)^{-1}||\leq C_1/\RE z$ 
for all $\RE z>a$. 
Thus the first statement follows from Lemma 3.5. 
Moreover, $\{x(xI-A_M)^{-1}\}_{x\geq a}$ converges to 0 in the strong 
topology as $x$ tends to $+\infty$. 
Thanks to Lemma 3.1 and Lemma 3.3, for $w\in \rH$ we have 
\begin{eqnarray*}
||x(xI-A_M)^{-1}e_w-x(xI-A)^{-1}e_w||&=&x|(\xi_{M,x},e_w)|||e_x||\\
&=&\sqrt{\frac{1}{2x}}\frac{|1-\frac{M(w)}{M(x)}|}{|1-\frac{w}{x}|} 
\to 0, \quad (x\to+\infty),
\end{eqnarray*}
which shows the statement. 
\end{proof}
\section{Global Theory}
Let $\{S_t\}_{t\geq 0}$ be as before and $\{K_t\}_{t\geq 0}$ be 
a family of bounded operators on $\Hi$. 
Then it is easy to show that the $C_0$-semigroup condition for 
$\{T_t=S_t+K_t\}$ is equivalent to (K0) and that the condition (C1) 
is equivalent to (K1) below respectively: 
\begin{itemize}
\item [(K0)] The map $\R_{>0}\ni t\mapsto K_t$ is strongly continuous and 
$$K_{s+t}=K_sS_t+S_sK_t+K_sK_t,\quad \forall s,t\geq 0,$$
\item [(K1)] 
$$K_t\Hi\subset L^2(0,t),\quad \forall t\geq 0.$$ 
\end{itemize}

We claim that the above conditions and (C2) altogether are equivalent to 
\begin{itemize}
\item [(K2)] There exists a measurable function $k(x,y)$ defined on 
$(0,\infty)^2$ satisfying 
$$K_tf(x)=\left\{
\begin{array}{ll}\int_0^tk(t-x,y)f(y)dy& (x<t)\\
0 &(t\leq x)
\end{array}\right.,$$
$$k(x+t,y)=k(x,y+t)+\int_0^tk(x,s)k(t-s,y)ds,\quad \forall t\geq 0, 
\;\textrm{a.e.}\; (x,y)\in (0,\infty)^2,$$
$$\int_0^t\int_0^\infty|k(x,y)|^2dydx<+\infty, \quad \forall t>0.$$
\end{itemize}

Indeed, assume that $\{T_t\}_{t\geq 0}$ is a $C_0$-semigroup satisfying 
(C1) and (C2). 
Let $k_t$ be the integral kernel for $K_t$. 
Then $k_t\in L^2(0,\infty)^2$ with support in 
$(0,t]\times (0,\infty)$. 
The semigroup property of $\{T_t\}$ implies that for every 
$s,t\geq 0$ and $f\in \Hi$, the following holds:  
\begin{eqnarray*}\lefteqn{\int_0^\infty k_{s+t}(x,y)f(x)dx}\\
&=&\int_t^\infty k_s(x,y)f(y-t)dy+1_{[s,s+t]}(x)
\int_0^\infty k_t(x-s,y)f(y)dy\\
&+&\int_0^\infty \int_0^t k_s(x,u)k_t(u,y)f(y)dudy\\
&=&\int_0^\infty \Big(k_s(x,y+t)+1_{[s,s+t]}(x) k_t(x-s,y)+
\int_0^t k_s(x,u)k_t(u,y)du\Big)f(y)dy,\\
\end{eqnarray*} 
where $1_I(x)$ denotes the characteristic function of an interval $I$. 
Therefore, for almost all $(x,y)\in (0,\infty)^2$, we get 
$$k_{s+t}(x,y)=k_s(x,y+t)+1_{[s,s+t]}(x) k_t(x-s,y)+
\int_0^t k_s(x,u)k_t(u,y)du.$$
We set $\tilde{k}_t(x,y)=k_t(t-x,y)$ for $x\in (0,t]$ and 
set $\tilde{k}_t(x,y)=0$ for $x>t$. 
Then the above relation is equivalent to 
\begin{eqnarray*}\tilde{k}_{s+t}(s+t-x,y)&=&
1_{(0,s]}(x)\Big(\tilde{k}_s(s-x,y+t)+
\int_0^t\tilde{k}_s(s-x,u)\tilde{k}_t(t-u,y)du\Big)\\
&+&1_{[s,s+t]}(x)\tilde{k}_t(s+t-x,y),
\end{eqnarray*}
for $0\leq x\leq s+t$. 
This shows that there exists a measurable function $k(x,y)$ on $(0,1)^2$ 
such that $\tilde{k}_t(x,y)=k(x,y)$ for $(x,y)\in (0,t]\times (0,\infty)$ 
and (K2) holds. 
The converse implication also follows from the same computation. 

\begin{lemma} Let $\{T_t\}_{t\geq 0}$ be a $C_0$-semigroup 
satisfying (C1) and (C2) and let $K_t$ and $k(x,y)$ be 
as above. 
We define an $\Hi$-valued function $\kappa(x)$ on $(0,\infty)$ by 
$\kappa(x)(\cdot)=k(x,\cdot)$.  
Then
\begin{itemize}
\item [$(1)$] There exists a positive number $a$ such that 
$e^{-ax}k(x,y)$ is square integrable on $(0,\infty)^2$. 
In particular, 
$$\Lt \kappa(z)(\cdot):=\int_0^\infty e^{-zx}k(x,\cdot)dx$$ 
makes sense as $\Hi$-valued holomorphic function on $\RE z>a$. 
\item [$(2)$] Let $A_M$ be the generator of $\{T_t\}_{t\geq 0}$ 
with $M\in \HD$. 
Then $\xi_{M,z}=\Lt \kappa(z)$ and 
$$\int_0^\infty \int_0^\infty |e^{-rx}k(x,y)|^2dxdy=
\frac{1}{2\pi}\int_{-\infty}^\infty||\xi_{M,r+is}||^2ds,\quad 
\forall r\geq a.$$ 
\end{itemize}
\end{lemma}

\begin{proof} (1) The condition (K0) implies 
$$||K_{s+t}||_{H.S.}\leq ||K_s||_{H.S.}+||K_t||_{H.S.}+
||K_s||_{H.S.}||K_t||,$$
where $||\cdot||_{H.S.}$ denotes the Hilbert-Schmidt norm. 
Thus for every natural number $n\in \N$, we get 
$$||K_{n+1}||_{H.S.}\leq ||K_n||_{H.S.}+||K_1||_{H.S.}+
||K_n||_{H.S.}||K_1||,$$which implies
$$||K_n||_{H.S.}\leq \frac{||K_1||_{H.S.}}{||K_1||}(1+||K_1||)^{n+1}.$$
Since 
$||K_t||_{H.S.}=\Big(\int_0^t\int_0^\infty|k(x,y)|^2dydx\Big)^{1/2}$ 
is an increasing function in $t$ with the above estimate for 
$t=n\in \N$, there exist two positive numbers $M$ and $b$ such 
that 
$$||K_t||_{H.S.}^2\leq Me^{bt}.$$
We choose $a$ satisfying $2a>b$ and set $\theta(t)=||K_t||_{H.S.}^2$. 
Then
\begin{eqnarray*}
\int_0^r\int_0^\infty |e^{-ax}k(x,y)|^2dydx&=&
\int_0^re^{-2at}\theta'(t)dt
=\theta(r)e^{-2ar}+2a\int_0^r\theta(t)e^{-2at}dt\\
&\leq&Me^{-(2a-b)r}+2a\int_0^rMe^{-(2a-b)t}dt\\
&\leq& \frac{2aM}{2a-b}, 
\end{eqnarray*}
which shows $e^{-ax}k(x,y)$ is square integrable. 

(2) We compute the Laplace transformation of $K_t$. 
For $f,g\in \Hi$ and $z$ with $\RE z>a$, we get
\begin{eqnarray*}\int_0^\infty e^{-zt}(K_tf,g)dt&=&
\int_0^\infty \int_0^t\int_0^\infty e^{-zt}k(t-x,y)f(y)g(x)
dydxdt\\
&=&\int_0^\infty \int_x^\infty\int_0^\infty e^{-zt}k(t-x,y)f(y)g(x)
dydtdx\\
&=&\int_0^\infty \int_0^\infty\int_0^\infty e^{-z(t+x)}k(t,y)f(y)g(x)
dydtdx\\
&=&(\Lt \kappa(z),f)(e_z,g). 
\end{eqnarray*}
Therefore Lemma 3.1 and Theorem 3.2 imply $\xi_{M,z}=\Lt \kappa(z)$. 
The last equation follows from the Plancherel theorem. 
\end{proof}

\begin{theorem} Let $M\in \HD$. 
Then $M\in \HD_2$ if and only if there exist positive numbers 
$a$ and $C$ such that $\xi_{M,z}$ is holomorphic on $\RE z> a$ 
and for all $x> a$
$$\int_{-\infty}^{+\infty}||\xi_{M,x+iy}||^2dy\leq C.$$
\end{theorem}

\begin{proof} Assume that $A_M$ generates a $C_0$-semigroup 
satisfying the conditions (C1) and (C2). 
Lemma 3.1 and Lemma 4.1 imply that there exist positive 
numbers $a$ and $M$ satisfying the above property. 

Assume now that there exist $a$ and $M$ satisfying the 
above property. 
Then Paley-Wiener theorem \cite[p. 163]{Y} implies that 
there exists a $\Hi$-valued measurable function $\kappa(x)$ on $(0,\infty)$ 
such that 
$$\xi_{M,z}=\int_0^\infty \kappa(x)e^{-zx}dx$$ and 
$e^{-ax}||\kappa(x)||$ is square integrable. 

We define a family of Hilbert-Schmidt operators $K_t\in \B(\Hi)$ for 
$t\geq 0$ by 
$$K_tf(x)=\left\{
\begin{array}{ll}(\kappa(t-x),f)& (x<t)\\
0 &(t\leq x)
\end{array}
\right.$$ 
Note that 
$$\int_0^\infty e^{-zt}K_tf=(f,\xi_{M,z})e_z$$
holds for all $f\in \Hi$ and $z$ with $\RE z>a$. 
It is easy to show that the map $t\mapsto K_t$ 
is continuous in the strong operator topology. 
Thus to show that $\{S_t+K_t\}_{t\geq 0}$ is a $C_0$-semigroup 
generated by $A_M$, it suffices to prove
$$K_{s+t}=K_sS_t+S_sK_t+K_sK_t, \quad \forall s,t\geq 0,$$
or equivalently to prove 
\begin{eqnarray*}\lefteqn{
\int_0^\infty ds\int_0^\infty dte^{-(z_1s+z_2t)}(K_{s+t}f,g)}\\
&=&\int_0^\infty ds\int_0^\infty dt e^{-(z_1s+z_2t)}
(K_sS_tf+S_sK_tf+K_sK_tf,g)
\end{eqnarray*}
for all $f,g$ in a total set in $\Hi$ and for all $\RE z_1>a$, $\RE z_2>a$. 
Assume $z_1\neq z_2$. 
Then the left-hand side is 
\begin{eqnarray*}\lefteqn{
\int_0^\infty ds\int_0^\infty dte^{-(z_1s+z_2t)}(K_{s+t}f,g)
=\int_0^\infty ds\int_s^\infty dte^{-(z_1s+z_2(t-s))}(K_tf,g)}\\
&=&\int_0^\infty dt(K_tf,g)e^{-z_2t}\int_0^t ds e^{-(z_1-z_2)s}
=\int_0^\infty dt(K_tf,g)\frac{e^{-z_2t}-e^{-z_1t}}{z_1-z_2}\\
&=&\frac{(f,\xi_{M,z_2})(e_{z_2},g)-(f,\xi_{M,z_1})(e_{z_1},g)}
{z_1-z_2}.
\end{eqnarray*}
The right-hand side is 

\begin{eqnarray*}\lefteqn{
((z_2I-A)^{-1}f,\xi_{M,z_1})(e_{z_1},g)
+(f,\xi_{M,z_2})((z_1I-A)^{-1}e_{z_2},g)}\\
&+&
(f,\xi_{M,z_2})(e_{z_2},\xi_{M,z_1})(e_{z_1},g)\\
&=&((z_2I-A)^{-1}f,\xi_{M,z_1})(e_{z_1},g)
+\frac{(f,\xi_{M,z_2})(e_{z_2},g)-(f,\xi_{M,z_2})
(e_{z_1},g)}{z_1-z_2}\\
&+&(f,\xi_{M,z_2})(e_{z_1},g)
\frac{M(z_1)-M(z_2)}{(z_1-z_2)M(z_1)}\\
&=&((z_2I-A)^{-1}f,\xi_{M,z_1})(e_{z_1},g)
+\frac{(f,\xi_{M,z_2})(e_{z_2},g)}{z_1-z_2}
-\frac{(f,\xi_{M,z_2})(e_{z_1},g)M(z_2)}{(z_1-z_2)M(z_1)},
\end{eqnarray*}
where we use $(\xi_{M,z},e_w)=\Lt{\xi_{M,z}}(w)$ and Lemma 3.3.
Thus, it suffices to show 
$$((z_2I-A)^{-1}f,\xi_{M,z_1})
=\frac{(f,\xi_{M,z_2})M(z_2)-(f,\xi_{M,z_1})M(z_1)}
{(z_1-z_2)M(z_1)}.$$ 
Setting $f=e_w$ with $\RE w>0$ and $w\neq z_1$, $w\neq z_2$, 
we can show this equality using Lemma 3.3, and so 
we get the statement. 
\end{proof}

\begin{cor} Let $M(z)$ be a holomorphic function on $\rH$.
\begin{itemize} 
\item [$(1)$] Assume that $M\in \HD_2$. then
\begin{itemize}
\item [$(\rm{i})$] 
$$\sup_{y\in \R}\frac{\Pin{|M(i\cdot)|^2}(x+iy)}{|M(x+iy)|^2}
=O(1),\quad (x\to+\infty),$$
\item [$(\rm{ii})$]
$$\int_{-\infty}^{+\infty}
\big(\frac{\Pin{|M(i\cdot)|^2}(x+iy)}{|M(x+iy)|^2}-1\big)dy=o(x),
\quad (x\to+\infty).$$
\item [$(\rm{iii})$]
There exists positive constant $a>0$ and $f\in L^2_{\mathrm{loc}}
(0,\infty)$ such that $e_af\in \Hi$ and $1/\big(zM(z)\big)$ is 
the Laplace transformation of $f$. 
In particular, there exists a positive constant $C$ such that 
$$|M(z)|\geq \frac{C\sqrt{\RE z-a}}{|z|},\quad \RE z>a.$$
\end{itemize}
\item [$(2)$] If $M$ satisfies $M(z)/(1+z)\in \Ha$ and  
$$\int_{-\infty}^{+\infty}
\big(\frac{\Pin{|M(i\cdot)|^2}(x+iy)}{|M(x+iy)|^2}-1\big)dy
=O(x),\quad (x\to+\infty),$$
then $M\in \HD_2$. 
\end{itemize}
\end{cor}

\begin{proof} (1) (i) follows from Corollary 3.6. 

Let the notation be as in the proof of the previous theorem. 
Then 
$$\frac{1}{2\pi}\int_{-\infty}^\infty||\xi_{M,x+iy}||^2dy
=\int_0^\infty||\kappa(s)||^2e^{-2xs}ds\rightarrow 0,
\quad (x\to+\infty),$$
which together with Lemma 3.5 and Lemma 4.1 implies (ii). 

Take $a>0$ such that 
$$\int_0^\infty\int_0^\infty|e^{-as}k(s,t)|^2dsdt<\infty.$$
For $x\geq a$, we define a Hilbert-Schmidt operator $L_x$ by 
$$L_xf(s)=\int_0^\infty e^{-xs}k(s,t)f(t)dy.$$
Lemma 3.3 and Lemma 4.1 imply
$$\Lt{L_xe_w}(z)=\frac{M(x+z)-M(w)}{(x+z-w)M(x+z)},$$
and so, 
$$||L_xe_w||^2=\frac{1}{2\pi}\int_{-\infty}^\infty
\big|\frac{1-\frac{M(w)}{M(x+iy)}}{x+iy-w}\big|^2dy.$$
On the other hand, 
$$||L_xe_w||\leq ||L_x||||e_w||\leq ||L_a||_{H.S.}||e_w||,$$ 
and we get 
$$\sup_{x\geq a}\int_{-\infty}^\infty\frac{dy}{|x+iy|^2|M(x+iy)|^2}
<\infty.$$
Thus (iii) follows from the Paley-Winer theorem. 

(2) Let $M$ be a function satisfying the assumption. 
If $M\notin \Ha$, we are done thanks to Theorem 4.2.  
Note that $\xi_{M,z}$ still makes sense for $M$ without 
the condition $M\notin \Ha$. 
Using the same argument as in the proof of Theorem 4.2, 
we can construct a $C_0$-semigroup $\{T_t\}_{t>0}$, whose 
Laplace transformation is given by 
$(zI-A)^{-1}+e_z\otimes \xi_{M,z}$. 
Let $A_{M_1}$ be the generator of the $\{T_t\}_{t>0}$. 
Then we have $\xi_{M,z}=\xi_{M_1,z}$, and so 
$$\frac{M(z)-M(w)}{M(z)(z-w)}=\frac{M_1(z)-M_1(w)}{M_1(z)(z-w)}$$
thanks to Lemma 3.3. 
Letting $w$ tend to $z$, we get $M'(z)/M(z)=M_1'(z)/M_1(z)$, 
and consequently $M$ is a scalar multiple of $M_1$, 
which shows $M\notin \Ha$. 
\end{proof}

\begin{cor} Let $M(z)$ be a holomorphic function on $\rH$ satisfying 
$M(z)/(1+z)\in \Ha$.  
We assume that there exist positive constants 
$a,n$, and $C$ such that 
$$\frac{1}{|M(z)|}\leq C(1+|z|)^n$$
holds for all $z$ with $\RE z\geq a$ and 
$$\int_{-\infty}^{+\infty}
\big(\frac{\Pin{|M(i\cdot)|^2}(x+iy)}{|M(x+iy)|^2}-1\big)dy<+\infty$$
for some $a<x$. 
Then $M\in \HD_2$. 
\end{cor}

\begin{proof} First we claim that there exists a distribution 
$k\in \cD(\R^2)'$ whose support is in $[0,\infty)^2$ such that 
the function $\Lt{\xi_{M,z}}(w)$ in $(z,w)$ is the Laplace transformation 
of $k(s,t)$.  
Indeed, let $q\in \Hi$ satisfying $(1+z)\Lt q(z)=M(z)$. 
Lemma 3.3 shows 
$$M(z)\Lt{\xi_{M,z}}(w)=(1+z)\frac{\Lt q(z)-\Lt q(w)}{z-w}+\Lt q(w).$$
Note that 
$$-\frac{\Lt q(z)-\Lt q(w)}{z-w}$$ 
is the Laplace transformation of the distribution $q(x+y)$. 
Therefore using the famous criterion \cite{S}, we can show 
the claim. 

By Lemma 3.5 and the assumption, the distribution $e^{-xs}k(s,t)$ in the two variables $(s,t)$ is actually 
a square integrable function, and so $k(s,t)$ belongs to $L^1_{\mathrm{loc}}[0,\infty)^2$. 
The rest of the proof is the same as that of Theorem 4.2. 
\end{proof}

Let $\D$ be the unit disc and $\zeta=\frac{z-1}{z+1}$. 
When $M\in \HD$, the function $\D\ni \zeta\mapsto M(z)$ is in the Hardy 
class $H^2(\D)$, and so $M$ uniquely factories in the form 
$M(z)=cB(z)S(z)F(z)$ with the following property \cite[p.132]{H}: 
$c$ is a complex number of modulus 1. 
$B(z)$ is the Blaschke product
$$B(z)=\Big(\frac{z-1}{z+1}\Big)^k\prod_{n}\frac{|1-\beta_n^2|}{1-\beta_n^2}
\cdot\frac{z-\beta_n}{z+\overline{\beta_n}},$$
where $k$ is a non-negative integer and $\{\beta_n\}$ is a 
(finite or infinite) sequence in $\rH$ satisfying 
$$\sum_n\frac{\RE \beta_n}{1+|\beta_n|^2}<+\infty.$$ 
The singular component $S(z)$ is given by 
$$S(z)=e^{-\sigma z}\exp\big[-\int_{-\infty}^\infty 
\frac{\lambda z+i}{\lambda+iz}d\mu(\lambda)\big],$$
where $\sigma$ is a non-negative real number and $\mu$ is a finite 
singular positive measure. 
$B(z)$ and $S(z)$ are inner functions in the sense that 
$|B(z)|<1$, $|S(z)|<1$ for all $z\in \rH$ and 
$|B(i\lambda)|=|S(i\lambda)|=1$ for almost every $\lambda$. 

The outer component $F(z)$ is given by 
$$F(z)=\exp\big[\frac{1}{2\pi}\int_{-\infty}^\infty 
\frac{\lambda z+i}{\lambda+iz}\cdot\frac{\rho(\lambda)
d\lambda}{1+\lambda^2}\big],$$
where $\rho(\lambda)=\log |M(i\lambda)|^2$. 

\begin{theorem} Let $M=cB(z)S(z)F(z)\in \HD$ be as above.
Then $M\in \HD_2$ if and only if $B,S,F\in \HD_2$. 
Moreover, 
\begin{itemize}
\item [$(1)$] $B\in \HD_2$ if and only if 
$$\sum_n\RE \beta_n<+\infty.$$
\item [$(2)$] $S\in \HD_2$ if and only if $\sigma=0$ and  
$$\int_{-\infty}^\infty (1+\lambda^2)d\mu(\lambda)<+\infty.$$
\item [$(3)$] $F\in \HD_2$ if and only if 
$$\sup_{y\in \R}\big(\log\Pin{e^\rho}(x+iy)-\Pin\rho(x+iy)\big)=O(1),
\quad (x\to+\infty),$$
$$\int_{-\infty}^{+\infty}\big(\log\Pin{e^\rho}(x+iy)
-\Pin\rho(x+iy)\big)dy=O(x),\quad (x\to+\infty).$$
\item [$(4)$] Assume that there exist positive constants $a$, $n$, and 
$C$ such that for all $\RE z\geq a$, 
$$\frac{1}{|F(z)|}\leq C(1+|z|)^n.$$
Then $F\in \HD_2$ if and only if there exists $x> a$ such that 
$$\sup_{y\in \R}\big(\log\Pin{e^\rho}(x+iy)-\Pin\rho(x+iy)\big)<+\infty,$$
$$\int_{-\infty}^{+\infty}\big(\log\Pin{e^\rho}(x+iy)
-\Pin\rho(x+iy)\big)dy<+\infty.$$
\end{itemize}
\end{theorem}

\begin{proof} Let 
$E_1(z)=-\log|B(z)|^2$, $E_2(z)=-\log|S(z)|^2$ and 
$$E(z)=\log\Pin{e^\rho}(z)-\Pin\rho(z)=\log\Pin{|F(i\cdot)|^2}(z)-
\log|F(z)|^2,$$
which are non-negative functions. 
Then we have 
\begin{eqnarray*}\lefteqn{\frac{\Pin{|M(i\cdot)|^2}(z)}{|M(z)|^2}-1
=e^{E_1(z)+E_2(z)+E_3(z)}-1}\\
&=& \big(e^{E_1(z)}-1\big)\big(e^{E_2(z)}-1\big)\big(e^{E_3(z)}-1\big)\\
&+&\big(e^{E_1(z)}-1\big)\big(e^{E_2(z)}-1\big)
+\big(e^{E_2(z)}-1\big)\big(e^{E_3(z)}-1\big)
+\big(e^{E_1(z)}-1\big)\big(e^{E_3(z)}-1\big)\\
&+&\big(e^{E_1(z)}-1\big)+\big(e^{E_2(z)}-1\big)+\big(e^{E_3(z)}-1\big)\\
&\geq& \big(e^{E_1(z)}-1\big)+\big(e^{E_2(z)}-1\big)+\big(e^{E_3(z)}-1\big)\\
&\geq& E_1(z)+E_2(z)+E_3(z).
\end{eqnarray*}
Thus Corollary 3.6 and Theorem 4.2 implies that 
$M\in \HD_2$ if and only if $B,S,F\in \HD_2$. 
(3) and (4) also follows from the above computation with an easy 
inequality 
$$x\leq e^x-1\leq \frac{e^r-1}{r}x,\quad 0\leq x\leq r. $$ 

(1) Clearly every inner function belongs to $\HD$. 
We have 
$$-\log|B(z)|^2=k\log|\frac{z+1}{z-1}|^2+\sum_n
\log|\frac{z+\overline{\beta_n}}{z-\beta_n}|^2.$$
Note that each term in the summation above is non-negative and integrable, 
and so 
$$\int_{-\infty}^\infty-\log|B(x+iy)|^2dy=k\int_{-\infty}^\infty
\log|\frac{x+1+iy}{x-1+iy}|^2dy+\sum_n \int_{-\infty}^\infty
\log|\frac{x+iy+\overline{\beta_n}}{x+iy-\beta_n}|^2dy.$$
Let $\beta_n=a_n+ib_n$ with $a_n,b_n\in \R$ and $a_n>0$. 
Then we have 
$$\int_{-\infty}^\infty
\log|\frac{x+iy+\overline{\beta_n}}{x+iy-\beta_n}|^2dy
=\pi\big(x+a_n-|x-a_n|).$$
The integral $\int_{-\infty}^\infty-\log|B(x+iy)|^2dy$ is finite if and only if $a_n<x$ except for 
finitely many $n$ and $\sum_na_n<\infty$. 
Thus if $B\in \HD_2$, we get $\sum_na_n<\infty$.  

Conversely, assume $\sum_na_n<\infty$. 
Let $a=\sup_na_n +2$ and $x\geq a$. 
Then 
\begin{eqnarray*}\log|\frac{x+iy+\overline{\beta_n}}{x+iy-\beta_n}|^2
&=&\log\frac{(x+a_n)^2+(y-b_n)^2}{(x-a_n)^2+(y-b_n)^2}\leq 
\frac{4xa_n}{(x-a_n)^2+(y-b_n)^2}\\
&\leq&\frac{4xa_n}{(x+1-a)^2},
\end{eqnarray*}
and so 
$$\frac{1}{|B(x+iy)|^2}\leq \exp\big[\frac{4x}{(x+1-a)^2}
\big(k+\sum_na_n\big)\big].$$
Thus there exists a positive constant $r$ such that 
$|B(x+iy)|^{-2}\leq r$. 
As we have 
$$\frac{1}{|B(x+iy)|^2}-1\leq -\frac{e^r-1}{r}\log|B(x+iy)|^2,$$
we get $B\in \HD_2$. 

(2) We have 
$$-\log|S(x+iy)|^2=2\sigma x+2\int_{-\infty}^\infty
\frac{x(1+\lambda^2)}{x^2+(y-\lambda)^2}d\mu(\lambda).$$ 
Since 
$$\int_{-\infty}^\infty\int_{-\infty}^\infty
\frac{x(1+\lambda^2)}{x^2+(y-\lambda)^2}d\mu(\lambda)dy
=\pi\int_{-\infty}^\infty(1+\lambda^2)d\mu(\lambda),$$
if $S\in \HD_2$, then $\sigma=0$ and $(1+\lambda^2)\mu$ is a 
finite measure. 

Conversely assume that $\sigma=0$ and $(1+\lambda^2)\mu$ is a 
finite measure. 
Then 
$$-\log|S(x+iy)|^2\leq\frac{2}{x}\int_{-\infty}^\infty
(1+\lambda^2)d\mu(\lambda),$$
and there exists a positive constant $l$ such that 
$|S(x+iy)|^{-2}\leq l$ for all $x\geq 1$. 
In the same way as above, this implies $S\in \HD_2$.   
\end{proof}

\begin{cor} If $M_1,M_2\in \HD_2$ are outer functions, then 
$M_1^\theta M_2^{1-\theta}\in \HD_2$ for all $0\leq \theta\leq 1$. 
\end{cor}

\begin{proof} 
H\"older inequality implies 
\begin{eqnarray*}\lefteqn{
\frac{\Pin{|M^\theta_1(i\cdot)M_2^{1-\theta}(i\cdot)|^2}(x+iy)}
{|M_1(x+iy)|^{2\theta}|M_1(x+iy)|^{2-2\theta}}-1}\\
&\leq&\Big(\frac{\Pin{|M_1(i\cdot)|^2}(x+iy)}
{|M_1(x+iy)|^2}\Big)^\theta
\Big(\frac{\Pin{|M_2(i\cdot)|^2}(x+iy)}
{|M_2(x+iy)|^2}\Big)^{1-\theta}-1\\
&\leq&\Big(\frac{\Pin{|M_1(i\cdot)|^2}(x+iy)}
{|M_1(x+iy)|^2}\Big)
\Big(\frac{\Pin{|M_2(i\cdot)|^2}(x+iy)}
{|M_2(x+iy)|^2}\Big)-1\\
&=&\Big(\frac{\Pin{|M_1(i\cdot)|^2}(x+iy)}
{|M_1(x+iy)|^2}-1\Big)\Big(\frac{\Pin{|M_2(i\cdot)|^2}(x+iy)}
{|M_2(x+iy)|^2}-1\Big)\\
&+&\Big(\frac{\Pin{|M_1(i\cdot)|^2}(x+iy)}
{|M_1(x+iy)|^2}-1\Big)+
\Big(\frac{\Pin{|M_2(i\cdot)|^2}(x+iy)}
{|M_2(x+iy)|^2}-1\Big).
\end{eqnarray*}
\end{proof}

We show, by example, that not every $M\in \HD$ is in $\HD_2$.  
For $\alpha\in \R$, we set $M(z)=(1+z)^\alpha$. 
Then $M\in \HD$ if and only if $-1/2\leq \alpha<1/2$. 
Corollary 3.6 shows that if $\alpha=-1/2$, the operator $A_M$ does not 
generate a $C_0$-semigroup at all.  

\begin{theorem} Let $-1/2 <\alpha<1/2$ and $M(z)=(1+z)^\alpha$. 
Then $A_M$ generates a $C_0$-semigroup satisfying 
the conditions (C1) and (C2) if and only if $\alpha=0$. 
\end{theorem}

\begin{proof}
Thanks to Lemma 3.5, we have 
\begin{eqnarray*}\lefteqn{
\int_{-\infty}^\infty
\Big[\frac{\Pin{|M(i\cdot)|^2}(x+iy)}{|M(x+iy)|^2}-1\Big]dy}\\
&=&\int_{-\infty}^\infty\Big[
\frac{1}{[(1+x)^2+y^2]^{\alpha}}
\frac{x}{\pi}
\int_{-\infty}^\infty\frac{(1+s^2)^{\alpha}}{x^2+(y-s)^2}ds-1
\Big]dy\\
&=&x\int_0^{+\infty}\Big[
\frac{1}{[(1+\frac{1}{x})^2+r^2]^\alpha\pi}
\int_{-\infty}^\infty
\frac{[\frac{1}{x^2}+(r+u)^2]^{\alpha}}{1+u^2}du-1\Big]dr.
\end{eqnarray*}
Let $f_x(r)$ be the integrand of the right-hand side, 
which is a non-negative function. 
If $M\in \HD_2$, we have 
$$\lim_{x\to +\infty}\frac{1}{x}\int_{-\infty}^\infty
\Big[\frac{\Pin{|M(i\cdot)|^2}(x+iy)}{|M(x+iy)|^2}-1\Big]dy=0.$$
On the other hand, 
\begin{eqnarray*}
\liminf_{x\to +\infty}\int_0^{+\infty}f_x(r)dr
&\geq& \int_0^{+\infty}\liminf_{x\to +\infty}f_x(r)dr\\
&=&\int_0^{+\infty}\Big[
\frac{1}{(1+r^2)^\alpha\pi}
\int_{-\infty}^\infty
\frac{|r+u|^{2\alpha}}{1+u^2}du-1\Big]dr.
\end{eqnarray*}
The integrand of the last term is a continuous non-negative 
function whose value at 0 is 
$$\frac{1}{\pi}
\int_{-\infty}^\infty
\frac{|u|^{2\alpha}}{1+u^2}du-1=
\frac{1}{\cos\big(\alpha\pi\big)}-1.$$
This is zero if and only if $\alpha=0$. 
Therefore for $\alpha\neq 0$, we get $M\notin \HD_2$. 
When $\alpha=0$, the operator $A_M$ is nothing but $A$, which generates 
$\{S_t\}_{t\geq 0}$. 
\end{proof}

\begin{remark} It is possible to show that 
$M(z)=(1+z)^\alpha\in \HD_b$ for $-1/2<\alpha<1/2$ 
(see Example 6.7 and Corollary 7.8). 
\end{remark}

\begin{theorem} 
Let $M(z)=\log^\alpha(a+z)$ with $\alpha>0$ and $a\geq 1$. 
Then $M\in \HD_2$. 
\end{theorem}

\begin{proof} 
Thanks to Corollary 4.4, it suffices to show
$$\int_{-\infty}^\infty\big[\frac{\Pin{|M(i\cdot)|^2}(x+iy)}{|M(x+iy)|^2}-1
\big]dy<\infty$$
for a fixed $x$. 
Since $M(z)$ is an outer function, we may and do assume that 
$\alpha$ is a natural number thanks to Corollary 4.6. 
To conclude the above integral is finite, 
it suffices to estimate $\Pin{|M(i\cdot)|^2}(x+iy)$ and 
$|M(x+iy)|^2$ when $y$ goes to $+\infty$. 
We assume that $y$ is a sufficiently large positive number from now. 

\begin{eqnarray*}\lefteqn{|M(x+iy)|^2=
\big[\log^2\sqrt{(a+x)^2+y^2}+\arctan^2\frac{y}{a+x}\big]^\alpha}\\
&=&\big[\big(\log y+\log\sqrt{1+\frac{(a+x)^2}{y^2}}\big)^2
+\big(\frac{\pi}{2}-\arctan\frac{a+x}{y}\big)^2\big]^\alpha\\
&=&\big[\log^2 y+\frac{\pi^2}{4}+O(\frac{1}{y})]^\alpha\\
&=&\big[\log^2y+\frac{\pi^2}{4}\big]^\alpha
\big[1+O(\frac{1}{y\log^2 y})\big]. 
\end{eqnarray*}
On the other hand 
\begin{eqnarray*}\lefteqn{\Pin{|M(i\cdot)|^2}(x+iy)
=\frac{x}{\pi}\int_{-\infty}^\infty\frac{
\big[\log^2\sqrt{a^2+(y+\lambda)^2}+\arctan^2(\lambda+y)\big]^\alpha}
{x^2+\lambda^2}d\lambda}\\
&\leq&\frac{x}{\pi}\int_{-\infty}^\infty\frac{
\big[\log^2\sqrt{a^2+(y+\lambda)^2}+\frac{\pi^2}{4}\big]^\alpha}
{x^2+\lambda^2}d\lambda\\
&=&\sum_{k=0}^\alpha
\left(\begin{array}{c}\alpha\\k\\ \end{array}\right)
\big(\frac{\pi^2}{4}\big)^{\alpha-k}
\frac{x}{\pi}\int_{-\infty}^\infty\frac{\log^{2k}\sqrt{a^2+(y+\lambda)^2}}
{x^2+\lambda^2}d\lambda\\
&=&\big[\log^2y+\frac{\pi^2}{4}\big]^\alpha
+\sum_{k=1}^\alpha
\left(\begin{array}{c}\alpha\\k\\ \end{array}\right)
\big(\frac{\pi^2}{4}\big)^{\alpha-k}
\frac{x}{\pi}\int_{-\infty}^\infty\frac{
\log^{2k}\sqrt{a^2+(y+\lambda)^2}-\log^{2k}y}{x^2+\lambda^2}d\lambda.
\end{eqnarray*}
Thus it suffices to show 
$$\int_{-\infty}^\infty\frac{
\log^n\sqrt{a^2+(y+\lambda)^2}-\log^n y}{x^2+\lambda^2}d\lambda
=O(\frac{\log^{n-2}y}{y}),\quad (y\to+\infty),$$
for every even integer $n\geq 2$. 
We set 
$$f_1(y)=\int_{-\infty}^\infty\frac{
\log^n\sqrt{a^2+(y+\lambda)^2}-\log^n|y+\lambda|}{x^2+\lambda^2}d\lambda,$$
$$f_2(y)=\int_{-\infty}^\infty\frac{
\log^n|y+\lambda|-\log^n y}{x^2+\lambda^2}d\lambda,$$
and estimate $f_1$ and $f_2$ separately. 

Let $g(\lambda)=\log^n\sqrt{a^2+(\lambda)^2}-\log^n|\lambda|$. 
Then $g(\lambda)=O(\log^{n-1}|\lambda|/|\lambda|^2)$ as $|\lambda|$ 
tends to infinity and $g\in L^1(\R)$. 
Fix a number $r$ satisfying $0<r<1$. 
Then
\begin{eqnarray*}|f_1(y)|&\leq&\int_{-\infty}^\infty
\frac{|g(\lambda+y)|}{x^2+\lambda^2}d\lambda
=\int_{-\infty}^{-ry}
\frac{|g(\lambda+y)|}{x^2+\lambda^2}d\lambda
+\int_{-ry}^\infty
\frac{|g(\lambda+y)|}{x^2+\lambda^2}d\lambda\\
&\leq&\frac{||g||_1}{x^2+r^2y^2}+\sup_{\lambda\geq (1-r)y}
\{|g(\lambda)|\}\int_{-ry}^\infty
\frac{1}{x^2+\lambda^2}d\lambda\\
&=&O(\frac{\log^{n-1}y}{y^2}),\quad (y\to+\infty). 
\end{eqnarray*}

For $f_2$ we have 
\begin{eqnarray*}f_2(y)&=&\int_{-\infty}^\infty
\frac{\big[\log y+\log|1+\frac{\lambda}{y}|\big]^n-\log^n y}
{x^2+\lambda^2}d\lambda\\
&=&\frac{1}{y}\int_{-\infty}^\infty
\frac{\big[\log y+\log|u|\big]^n-\log^n y}
{\frac{x^2}{y^2}+(u-1)^2}du\\
&=&\sum_{k=1}^n\left(\begin{array}{c}n\\k\\ \end{array}\right)
\frac{\log^{n-k}y}{y}\int_{-\infty}^\infty\frac{\log^k|u|}
{\frac{x^2}{y^2}+(u-1)^2}du\\
&=&\frac{n\log^{n-1}y}{y}\int_{-\infty}^\infty\frac{\log|u|}
{\frac{x^2}{y^2}+(u-1)^2}du+O(\frac{\log^{n-2}y}{y}).\\
\end{eqnarray*}
Let $g_1(y)$ be the integral in the last term. 
To finish the proof, it suffices to show 
$g_1(y)=O(1/\log y)$. 
Indeed, 
\begin{eqnarray*}|g_1(y)|
&=&
\big|\int_{-1}^1\frac{\log|u|}{\frac{x^2}{y^2}+(u-1)^2}du
+\int_{|u|\geq 1}\frac{\log|u|}
{\frac{x^2}{y^2}+(u-1)^2}du\big|\\
&=&\big|\int_{-1}^1\frac{\log|u|}{\frac{x^2}{y^2}+(u-1)^2}du
-\int_{-1}^1\frac{\log|u|}
{\frac{x^2u^2}{y^2}+(u-1)^2}du\big|\\
&\leq&\frac{x^2}{y^2}\int_{-1}^1
\frac{(1+u)|\log|u||}{(1-u)[\frac{x^2}{y^2}+(1-u)^2]}du\\
&\leq&\frac{x^2}{y^2}\int_{-1}^{1/2}
\frac{(1+u)|\log|u||}{(1-u)^3}du
+\frac{x^2}{y^2}\int_{1/2}^{1-1/y}
\frac{(1+u)|\log|u||}{(1-u)^3}du\\
&+&\int_{1-1/y}^1
\frac{(1+u)|\log|u||}{(1-u)}du\\
&\leq&\frac{x^2}{y^2}\int_{-1}^{1/2}
\frac{(1+u)|\log|u||}{(1-u)^3}du+\frac{Cx^2(y-2)}{y^2}
+\frac{C}{y}=O(\frac{1}{y}),
\end{eqnarray*}
where 
$$C=\sup_{1/2<u<1}\frac{(1+u)|\log|u||}{(1-u)}.$$ 
This finishes the proof. 
\end{proof}

We end this section by stating the following fact essentially proven 
in Lemma 4.1. 
\begin{cor}
Let $M\in \HD_2$ with $e^{tA_M}=S_t+K_t$. 
Then 
$$\int_0^\infty e^{-2tx}d_t||K_t||_{H.S.}^2=
\frac{1}{2\pi}\int_{-\infty}^{+\infty}||\xi_{M,x+iy}||^2dy.$$
\end{cor}

It is well known that the asymptotic behavior of the above quantity for 
large $x$ determines that of $||K_t||_{H.S.}^2$ for small $t$. 
As we suspect that the latter might survive as a characteristic property 
of the $E_0$-semigroup arising from $e^{tA_M}$ (see Remark 5.8), 
it might be an interesting problem to determine the former for concrete 
examples. 
We will do it for a class of examples in next section. 

\section{A class of examples}
To introduce a class of functions in $\HD_2$, we give a heuristic 
argument first. 
In the sequel, for $\varphi,\psi \in L^1_{\mathrm{loc}}[0,\infty)$ we denote by $\varphi*\psi$ the convolution 
$$\varphi*\psi(x)=\int_0^x\varphi(s)\psi(x-s)ds=\int_0^x\varphi(x-s)\psi(s)ds.$$
Then we have $\Lt{\varphi*\psi}(z)=\Lt \varphi(z)\Lt \psi(z)$ whenever the both sides make sense. 

Let $k$ be a function satisfying the condition (K2). 
We set $l(x,z)$ to be the Laplace transform of $k(x,y)$ in the second 
variable, that is $l(x,z)=(\kappa(x),e_z)$. 
Then the condition (K2) implies 
$$l(x+t,z)=e^{zt}l(x,z)-e^{zt}\int_0^tk(x,y)e^{-zy}dy
+\int_0^t k(x,s)l(t-s,z)ds.$$
Assume for the moment that $k(x,y)$ is in $C^1$-class.  
Differentiating the both sides of the above equation by $t$ and 
setting $t=0$, we get 
$$\frac{\partial l}{\partial x}(x,z)=zl(x,z)+k(x,0)(l(0,z)-1).$$  
Solving this equation, we obtain 
$$e^{-zx}l(x,z)=l(0,z)+(l(0,z)-1)\int_0^xe^{-zs}k(s,0)ds,$$
which implies the following via the inverse Laplace transform: 
$$1_{(x,\infty)}(y)k(x,y-x)=\varphi(y)+(\varphi*1_{[0,x)}\psi)(y)-1_{[0,x)}(y)\psi(y),$$
where $\varphi(x)=k(0,x)$ and $\psi(x)=k(x,0)$. 
This is equivalent to 
$$k(x,y)=\varphi(x+y)+\int_0^x\varphi(y+x-s)\psi(s)ds,$$
and $\psi=\varphi+\varphi*\psi$. 

Now we start with $\varphi$ in an appropriate function space and 
we {\it define} $k$ by the above relation, which will turn out to 
satisfy the condition (K2). 
For two real numbers $x$ and $y$, we denote $x\wedge y=\min\{x,y\}$. 

\begin{lemma} Let $\varphi\in L^1_{\mathrm{loc}}[0,\infty)\cap 
L^2((0,\infty),(1\wedge x)dx)$. 
Then there exist a unique function 
$\psi\in L^1_{\mathrm{loc}}[0,\infty)$ 
and a positive number $a$ such that 
$e_a\psi\in L^1(0,\infty)$ and 
$$\Lt\psi(z)=\frac{\Lt\varphi(z)}{1-\Lt\varphi(z)}$$
for all $z$ with $\RE z\geq a$. 
\end{lemma}

\begin{proof} Take $\varepsilon>0$ such that 
$$\int_0^\varepsilon |\varphi(x)|dx<1.$$
Then, for $a>0$ we have 
$$||e_a\varphi||_1\leq \int_0^\varepsilon |\varphi(x)|dx 
+\int_\varepsilon^\infty |\varphi(x)|e^{-ax}dx
\leq \int_0^\varepsilon |\varphi(x)|dx+
\sqrt{\frac{e^{-2\varepsilon a}}{2a}
\int_\varepsilon^\infty |\varphi(x)|^2dx}.$$
Thus we can choose large $a$ so that $||e_a\varphi||_1<1.$
Note that the $n$-fold convolution of $e_a\varphi$ is 
$e_a\varphi^{*n}$. 
Thus, 
$$\sum_{n=1}^\infty e_a\varphi^{*n}$$
converges in $L^1(0,\infty)$ or equivalently 
$$\sum_{n=1}^\infty\varphi^{*n}$$
converges in $L^1((0,\infty),e^{-ax}dx)$. 
We define $\psi$ to be this limit. 
Then for $\RE z\geq a$ we get 
$$\int_0^\infty \psi(x)e^{-zx}dx=\sum_{n=1}^\infty\Lt\varphi(z)^n
=\frac{\Lt\varphi(z)}{1-\Lt\varphi(z)}.$$
\end{proof}

We fix $\varphi$ and $\psi$ as above and set 
$\Phi(z)=\Lt\varphi(z)$ and $\Psi(z)=\Lt\psi(z)$. 
We assume that $\varphi$ is a non-zero function. 
Since $\Phi-\Psi+\Phi\Psi=0$ holds, we have 
$\varphi-\psi+\varphi*\psi=0$. 

\begin{theorem} Let $\varphi\in L^1_{\mathrm{loc}}[0,\infty)\cap 
L^2((0,\infty),(1\wedge x)dx)$ and let $\psi$ be the function determined by Lemma 5.1. 
We set 
$$k(x,y)=\varphi(x+y)+\int_0^x\varphi(x+y-s)\psi(s)ds.$$
Then $k(x,y)$ satisfies the condition (K2).
\end{theorem}
\begin{proof} Note that $k(x,y)$ makes sense as an element of 
$L^1_{\mathrm{loc}}[0,\infty)^2$. 
Let $a\geq 1$ be a positive number such that 
$e_a\psi\in L^1(0,\infty)$. 

First we show that $e^{-ax}k(x,y)$ is in $L^2(0,\infty)^2$. 
Indeed for the first term, we have
\begin{eqnarray*}
\int_0^\infty dx\int_0^\infty dy|e^{-ax}\varphi(x+y)|^2
&=&\int_0^\infty e^{-2ax}dx\int_x^\infty dy|\varphi(y)|^2\\
&=&\int_0^\infty |\varphi(y)|^2dy\int_0^y dxe^{-2ax}\\
&=&\frac{1}{2a}\int_0^\infty (1-e^{-2ay})|\varphi(y)|^2dy\\
&\leq&\int_0^\infty|\varphi(y)|^2(1\wedge y)dy<+\infty, 
\end{eqnarray*}
where we use the fact $(1-e^{-by})/b\leq 1\wedge y$ for all 
$y\geq 0$ and $b\geq 1$. 

Using $\int_0^x|\psi(s)|ds\leq e^{ax}||e_a\psi||_1$, 
we get the following estimate for the second term: 
\begin{eqnarray*}\lefteqn{
\int_0^\infty dx\int_0^\infty dy|e^{-ax}
\int_0^x\psi(s)\varphi(x+y-s)ds|^2}\\
&=&\int_0^\infty e^{-2ax}dx\int_0^\infty dy
\int_0^x\psi(r)\varphi(x+y-r)dr\int_0^x
\overline{\psi(s)\varphi(x+y-s)}ds\\
&\leq&  \int_0^\infty e^{-2ax}dx
\int_0^x|\psi(r)|dr\int_0^x|\psi(s)|ds
\int_0^\infty |\varphi(x+y-r)\varphi(x+y-s)|dy\\
&\leq&\int_0^\infty e^{-2ax}dx\int_0^x|\psi(r)|dr\int_0^x|\psi(s)|ds
\int_0^\infty \frac{|\varphi(x+y-r)|^2+|\varphi(x+y-s)|^2}{2}dy\\
&=&\int_0^\infty e^{-2ax}dx\int_0^x|\psi(r)|dr\int_0^x|\psi(s)|ds
\int_0^\infty |\varphi(x+y-s)|^2dy\\
&=&  ||e_a\psi||_1\int_0^\infty e^{-ax}dx
\int_0^x|\psi(x-s)|ds\int_0^\infty |\varphi(s+y)|^2dy\\
&=&  ||e_a\psi||_1\int_0^\infty ds \int_s^\infty dx|\psi(x-s)|e^{-ax}
\int_0^\infty |\varphi(s+y)|^2dy\\
&=&  ||e_a\psi||_1\int_0^\infty e^{-as}ds \int_0^\infty dx
|\psi(x)|e^{-ax}\int_0^\infty |\varphi(s+y)|^2dy\\
&=&  ||e_a\psi||_1^2\int_0^\infty e^{-as}ds 
\int_0^\infty |\varphi(s+y)|^2dy\\
&=&  ||e_a\psi||_1^2\frac{1}{a}
\int_0^\infty(1-e^{-ay}) |\varphi(y)|^2dy\\
&\leq&  ||e_a\psi||_1^2
\int_0^\infty|\varphi(y)|^2(1\wedge y)dy<+\infty.\\
\end{eqnarray*}
Thus we get
$$\int_0^\infty\int_0^\infty |e^{-ax}k(x,y)|^2dxdy\leq (1+||e_a\psi||_1)^2
\int_0^\infty|\varphi(x)|^2(1\wedge x)dx.$$

As before, we define an $\Hi$-valued function $\kappa(x)$ by 
$\kappa(x)(\cdot)=k(x,\cdot)$ and define $K_t\in \B(\Hi)$ by 
$$K_tf(x)=\left\{
\begin{array}{ll}(\kappa(t-x),f)& (x<t)\\
0 &(t\leq x)
\end{array}
\right.$$ 
Then we get  
\begin{eqnarray*}||K_t||^2&\leq& ||K_t||_{H.S.}^2
=\int_0^tdx\int_0^\infty dy|k(x,y)|^2\leq 
e^{2at}\int_0^\infty dx\int_0^\infty dy|e^{-ax}k(x,y)|^2\\
&\leq& e^{2at}(1+||e_a\psi||_1)^2
\int_0^\infty|\varphi(y)|^2(1\wedge y)dy.
\end{eqnarray*}
Thus for $\RE z> a$ and $f\in \Hi$, the Laplace transformation 
of $\{K_t\}_{t\geq 0}$ is well-defined and  
$$\int_0^\infty e^{-zt}K_tfdt
=\int_0^\infty e^{-zt}(\kappa(t),f)dt e_z.$$
Using the definition of $\kappa$, we get 
\begin{eqnarray*}\lefteqn{
\int_0^\infty e^{-zt}(\kappa(t),f)dt}\\
&=&\int_0^\infty dt e^{-zt}
\Big(\int_0^\infty dy\varphi(t+y)f(y)+
\int_0^\infty dy\int_0^t ds\psi(s)\varphi(y+t-s)f(y)\Big)\\
&=&\int_0^\infty f(y)dy\Big(\int_0^\infty dte^{-zt}\varphi(t+y)+
\int_0^\infty \psi(s)ds\int_s^\infty dt e^{-zt}\varphi(y+t-s)
\Big)\\
&=&\int_0^\infty f(y)dy\Big(\int_0^\infty dte^{-zt}\varphi(t+y)+
\int_0^\infty \psi(s)e^{-zs}ds\int_0^\infty dt e^{-zt}\varphi(y+t)
\Big)\\
&=&(1+\Psi(z))\int_0^\infty f(y)dy
\int_0^\infty dte^{-zt}\varphi(t+y)\\
&=&(1+\Psi(z))(\eta_{\varphi,z},f),
\end{eqnarray*}
where 
$$\eta_{\varphi,z}(x)=\int_0^\infty e^{-zt}\varphi(t+x)dt.$$
As before we get 
$$\Lt{\eta_{\varphi,z}}(w)=\frac{\Phi(w)-\Phi(z)}{z-w}.$$
In the same way as in the proof of Theorem 4.2, 
we can show
$$\int_0^\infty ds\int_0^\infty dt e^{-sz_1-tz_2}K_{s+t}e_w
=\int_0^\infty ds\int_0^\infty dt e^{-sz_1-tz_2}
(K_sK_t+K_sS_t+S_sK_t)e_w,$$
which finishes the proof. 
\end{proof}

Thanks to Theorem 3.2 there exists $M\in \HD$ such that 
$A_M$ is the generator of the semigroup 
$\{T_t=S_t+K_t\}_{t\geq 0}$ constructed above. 
The above argument shows that we have the relation
$$(1+\Psi(z))\eta_{\varphi,z}(x)=\xi_{M,z}(x).$$
Let $q\in \Hi\setminus D(A)$ such that $(1+z)\Lt q(z)=M(z)$. 
Since the left-hand side is continuous in $x$, so is the right-hand 
side and Equation (3.2) implies that $q$ is continuous. 
We set $x=0$. 
Then the left hand side is 
$(1+\Psi(z))\Phi(z)=\Psi(z)$.
On the other hand the right-hand side is $q(0)/M(z)-1$, 
which shows $1+\Psi(z)=q(0)/M(z)$. 
Since $\Psi$ is not a constant function, $q(0)\neq 0$.  
Now the equation $\Phi-\Psi+\Phi\Psi=0$ implies 
$M(z)=q(0)(1-\Phi(z))$. 

Summing up the argument as above, we get
\begin{theorem} Let $\varphi\in 
L^1_{\mathrm{loc}}[0,\infty)\cap L^2((0,\infty),(1\wedge x)dx)$ and 
let $\{T_t\}_{t\geq 0}$ be the $C_0$-semigroup 
constructed from $\varphi$ by the above argument. 
Then $\{T_t\}_{t\geq 0}$ is a $C_0$-semigroup satisfying 
the conditions (C1) and (C2). 
Let $A_M\in \HD_2$ be the generator of $\{T_t\}_{t\geq 0}$, and 
let $q\in \Hi\setminus D(A)$ be the function satisfying 
$(1+z)\Lt q(z)=M(z)$. 
Then $q$ is continuous at 0 and $M(z)=q(0)(1-\Lt\varphi(z))$. 
The relation between $\varphi$ and $q$ is given by  
$$q(x)=q(0)\Big(e^{-x}-\int_0^x \varphi(s)e^{s-x}ds\Big),$$ 
$$\varphi=-\frac{q+q'}{q(0)}.$$
\end{theorem}

\begin{remark} To construct an example of a $C_0$-semigroup 
satisfying only (C1), the above argument still works for 
the Dirac function $\varphi(x)=r\delta_a(x)$ with $r\in \C$ and 
$a>0$ if each argument is appropriately interpreted.  
Indeed, it is possible to show $M(z)=1-\Phi(z)=1-re^{-az}$ 
belongs to $\HD_b$ (Corollary 7.7). 
\end{remark}

\begin{example}
When $\varphi(x)=ce^{-dx}$ with $c\neq 0$ and $\RE d>0$, 
by easy computation we have $\Phi(z)=c/(z+d)$, $\Psi(z)=c/(z+d-c)$, 
$\psi(x)=ce^{(c-d)x}$, and $k(x,y)=ce^{(c-d)x-dy}$. 
\end{example}

In the rest of this section, we give an estimate of 
$||K_t||_{H.S.}$ for small $t$ for the class of examples we have 
constructed in this section. 
 
\begin{lemma} Let $\varphi\in L^1_{\mathrm{loc}}[0,\infty)\cap 
L^2((0,\infty),(1\wedge x)dx)$. 
We set $\Phi=\Lt\varphi$, $M=1-\Phi\in \HD_2$, and $K_t=e^{tA_M}-S_t$. 
Then 
$$\frac{1}{2\pi}\int_{-\infty}^{+\infty}||\xi_{M,x+iy}||^2dy
\sim\frac{1}{2x}\int_0^\infty (1-e^{-2xs})|\varphi(s)|^2ds,
\quad (x\to+\infty).$$
\end{lemma}

\begin{proof} Since $\varphi$ belongs to $L^1(0,\infty)+\Hi$, 
the quantity $\Phi(x+iy)$ converges to 0, uniformly in $y$, 
as $x$ tends to $+\infty$. 
Thus Lemma 3.3 implies 
\begin{eqnarray*}\lefteqn{
\frac{1}{2\pi}\int_{-\infty}^{+\infty}||\xi_{M,x+iy}||^2dy}\\
&=&\frac{1}{4\pi^2}\int_{-\infty}^\infty\frac{1}{|1-\Phi(x+iy)|^2}
\int_{-\infty}^\infty 
\frac{|\Phi(i\lambda)-\Phi(x+iy)|^2}{|x+iy-i\lambda|^2}d\lambda dy\\
&\sim&\frac{1}{4\pi^2}\int_{-\infty}^\infty\int_{-\infty}^\infty
\frac{|\Phi(i\lambda)-\Phi(x+iy)|^2}{|x+iy-i\lambda|^2}d\lambda dy\\
&=&\frac{1}{4\pi^2 x}\int_{-\infty}^\infty\int_{-\infty}^\infty
\frac{|\Phi(ixu+iy)-\Phi(x+iy)|^2}{1+u^2}dydu.
\end{eqnarray*}
Let $f(s)=\varphi(s)(e^{-isxu}-e^{-sx})$. 
Then $f\in \Hi$ and 
$$\Phi(ixu+iy)-\Phi(x+iy)=\int_0^\infty f(s)e^{-isy}ds.$$
Thus the Plancherel theorem implies that the last term above  
equals to 
\begin{eqnarray*}\lefteqn{
\frac{1}{2\pi x}\int_{-\infty}^\infty\int_{-\infty}^\infty
\frac{|f(s)|^2}{1+u^2}dsdu}\\
&=&\frac{1}{2\pi x}\int_{-\infty}^\infty\int_{-\infty}^\infty
\frac{|\varphi(s)|^2\big(1+e^{-2sx}-2e^{-sx}\cos(sxu)\big)}
{1+u^2}dsdu\\
&=&\frac{1}{2x}\int_0^\infty(1-e^{-2sx})|\varphi(s)|^2ds.
\end{eqnarray*}
\end{proof}

\begin{cor} Assume $\varphi\in \Hi$ and $M=1-\Lt\varphi$. 
Then, 
$$\frac{1}{2\pi}\int_{-\infty}^{+\infty}||\xi_{M,x+iy}||^2dy
\sim \frac{||\varphi||^2}{2x}, \quad (x\to+\infty),$$
$$||K_t||_{H.S.}^2\sim ||\varphi||^2 t,\quad (t\to+0).$$
\end{cor}

\begin{proof} 
The first statement immediately follows from Lemma 5.5, which together with 
Corollary 4.10 implies the second one via the Tauberian theorem 
\cite[Chapter XIII.5, Theorem 3]{F}. 
\end{proof}

\begin{remark}
The above argument actually shows that $\varphi\in \Hi$ if and only if 
$||K_t||_{H.S.}^2=O(t)$, $(t\to+0)$. 
On the other hand, we show in the subsequent paper \cite{IS} that 
the $E_0$-semigroup arising from $e^{tA_M}$ is of type I if and only if 
$\varphi\in \Hi$. 
In view of this fact, it is tempting to conjecture that there exists a 
cocycle conjugacy invariant of $E_0$-semigroups that is computable from 
the asymptotic behavior of $||K_t||_{H.S.}^2$ for small $t$. 
\end{remark}

Now we treat $\varphi\in L^1_{\mathrm{loc}}[0,\infty)\cap 
L^2((0,\infty),(1\wedge x)dx)\setminus \Hi$.
We set 
$$m_\varphi(t)=\int_t^\infty|\varphi(x)|^2dx,$$
which diverges as $t$ tends to 0. 
Note that $m_\varphi(t)$ is integrable on every finite interval of the 
form $(0,a)$, $a>0$ because the Fubini theorem implies 
$$\int_0^am_\varphi(s)ds
=\int_0^\infty (u\wedge a)|\varphi(u)|^2du
=\int_0^as|\varphi(s)|^2ds+a\int_a^\infty |\varphi(s)|^2ds.$$

Recall that a measurable function $L(t)$ on $(0,\infty)$ is said to be 
slowly varying at 0 if $L(t)$ is positive on a non-empty interval $(0,a)$ 
and for every $s>0$,
$$\lim_{t\to+0}\frac{L(st)}{L(t)}=1.$$
A measurable function $f(t)$ on $(0,\infty)$ is said to be 
regularly varying at 0 if there exist a real number $\alpha$ and 
a slowly varying function $L(t)$ such that $f(t)=t^\alpha L(t)$. 

\begin{lemma}Let the notation be as in Lemma 5.6. 
Then 
$$\int_0^\infty e^{-xs}d_s||K_s||_{H.S.}^2\sim 
\int_0^\infty e^{-xs}m_\varphi(s)ds, \quad (x\to+\infty).$$
In particular, if the function $t\mapsto \int_0^tm_\varphi(s)ds$ is  
regularly varying at 0,  
$$||K_t||_{H.S.}^2\sim \int_0^tm_\varphi(u)du
=\int_0^ts|\varphi(s)|^2ds+t\int_t^\infty 
|\varphi(s)|^2ds,\quad (t\to+0).$$
\end{lemma}

\begin{proof} First we claim that $tm_\varphi(t)$ converges to 0 as 
$t$ tends to 0. 
Indeed, fix a positive number $\varepsilon$. 
Then for $0<t\leq\ \varepsilon$
$$tm_\varphi(t)\leq t\int_{\varepsilon}^\infty|\varphi(s)|^2ds
+\int_t^\varepsilon s|\varphi(s)|^2ds,$$
and 
$$\limsup_{t\to+0}tm_\varphi(t)\leq \int_0^\varepsilon s|\varphi(s)|^2ds,$$
which shows the claim. 
Thanks to Corollary 4.10, Lemma 5.5, and the above claim, 
the first statement follows from integration by part. 
The second statement follows from \cite[Chapter XIII.5, Theorem 3]{F}. 
\end{proof}

It is a routine work to show the following from Lemma 5.9 and 
\cite[Chapter VIII.9]{F}. 

\begin{cor} Let the notation be as in Lemma 5.6 and let 
$|\varphi(x)|=x^{\alpha-1}L(x)$ where $\alpha$ is a non-zero 
constant with $0\leq \alpha\leq 1/2$ and $L(x)$ is a slowly varying 
function at 0. 
\begin{itemize}
\item [$(1)$] When $\alpha=0$, 
$$||K_t||_{H.S.}^2\sim \int_0^t \frac{L(s)^2}{s}ds,\quad(t\to+0).$$
In particular, when $|\varphi(x)|\sim Cx^{-1}(\log x^{-1})^{-\beta}$, 
$(x\to+0)$ with $\beta>1$ and a non-zero constant $C$, 
$$||K_t||_{H.S.}^2\sim \frac{|C|^2}
{(2\beta-1)(\log\frac{1}{t}\big)^{2\beta-1}},\quad(t\to+0).$$
\item [$(2)$] When $0<\alpha<1/2$, 
$$||K_t||_{H.S.}^2\sim \frac{t^{2\alpha}L(t)^2}
{2\alpha(1-2\alpha)},\quad (t\to+0).$$
\item [$(3)$] When $\alpha=1/2$, 
$$||K_t||_{H.S.}^2\sim t\int_t^\infty\frac{L(s)^2}{s}ds,\quad (t\to+0).$$ 
In particular, when $|\varphi(x)|\sim Cx^{-1/2}(\log x^{-1})^\beta$, 
$(x\to +0)$ with $\beta\geq -1/2$ and a non-zero constant $C$, 
$$||K_t||_{H.S.}^2\sim \left\{
\begin{array}{ll}
\frac{|C|^2 t\big(\log\frac{1}{t}\big)^{2\beta+1}}{2\beta+1}
&(\beta>-\frac{1}{2})\\
|C|^2t\log\log\frac{1}{t}& (\beta=-\frac{1}{2})
\end{array}
\right.,\quad (t\to+0).$$
\end{itemize} 
\end{cor}

\section{Off white noise models}
The purpose of this section is to clarify the relationship between 
Tsirelson's off white noise and our semigroup $\{T_t\}_{t\geq 0}$. 

Let $M\in \HD$ be an outer function satisfying 
$\overline{M(\overline{z})}=M(z)$. 
We introduce a measure $\nu$ on $\R$ by setting 
$$\nu(d\lambda)=\frac{|M(i\lambda)|^2}{2\pi}d\lambda.$$ 
Let $G:=L^2(\R,\nu)$. 
We denote by $\inpr{\cdot}{\cdot}_G$ and $||\cdot||_G$ 
the inner product and the norm of $G$. 
Let $F$ be the unitary from $L^2(\R)$ to $G$ given by 
$Ff(\lambda)=\hat{f}(\lambda)/M(i\lambda)$, where 
$\hat{f}$ is the Fourier transform of $f$ 
$$\hat{f}(\lambda)=\Ft{f}(\lambda)
=\int_{-\infty}^\infty f(s)e^{-is\lambda}ds.$$
We regard $\Hi$ as a subspace of $L^2(\R)$ and so 
$Ff(\lambda)=\Lt f(i\lambda)/M(i\lambda)$ for $f\in \Hi$. 
We denote by $U_t$ the bilateral shift on $L^2(\R)$, that is, 
$U_tf(x)=f(x-t)$ and set $V_t=FU_tF^{-1}$. 
Then $V_t$ is given by the multiplication operator of $e^{-i\lambda t}$. 

Let $\cI_0$ be the set of finite open intervals of $\R$ and 
let $\cI$ be the set of (finite or infinite) open intervals of $\R$. 
We denote by $\cS=\cS(\R)$ the set of rapidly decreasing $C^\infty$ 
functions on $\R$. 
For $I\in \cI$, we set $\cS_I$ to be the linear span of 
$$\{1_Jf;\; f\in \cS,\; J\subset I,\; J\in \cI\}.$$

\begin{lemma} Let $I \in\cI$ and $f\in \cS$. 
Then $\Ft{1_If}\in G$. 
\end{lemma}

\begin{proof} Let $f\in \cS$ and $I$ be a finite interval $(a,b)$. 
Then 
$$\cF[1_If](\lambda)=\frac{\cF[1_If'](\lambda)+f(a)e^{-ia\lambda}
-f(b)e^{-ib\lambda}}{i\lambda}.$$
Since $M(z)/(1+z)\in \Ha$, we conclude that 
$\cF[1_If](\lambda)M(i\lambda)$ is square integrable. 
The case where $I$ is an infinite interval can be treated in the same way. 
\end{proof}

For $I\in \cI$, we set $G_0(I)=\Ft{\cS_I}$    
and set $G(I)$ to be the closure of $G_0(I)$. 
By definition, we have $V_tG(I)=G(t+I)$. 
We denote $G_0(-\infty,\infty)$ by $G_0$.  

\begin{lemma} With the notation as above, we have 
$G=G(-\infty,\infty)$ and $F\Hi=G(0,\infty)$.  
\end{lemma}

\begin{proof} Let $f\in G\cap G(-\infty,\infty)^\perp$. 
Then $fM$ is orthogonal to $gM(i\lambda)$ in $L^2(i\R)$ 
for all $g\in \cS$. 
Let $H_0$ be the closed linear span of $\{gM\}_{g\in \cS}$ 
in $L^2(i\R)$. 
Then $H_0$ is invariant under multiplication by $e^{-it\lambda}$ for 
all $t$ and so there exists a measurable subset $E\subset i\R$ such 
that $H_0=L^2(E)$. 
Since $M(z)/(1+z)\in \Ha$, the set of zeros of $M$ 
on the imaginary axis has Lebesgue measure zero and we can conclude 
$H_0=L^2(\R)$. 
Thus we get $fM=0$ and so $f=0$. 
The case of the half line can be treated in a similar way by using 
the fact that $M$ is an outer function and the Beurling-Lax theorem 
\cite[p.107]{H}. 
\end{proof}

We introduce an operator $Q_0$, called the {\it Riesz projection}, as follows: 
the domain $D(Q_0)$ is $G_0$ and 
$$Q_0\Ft{1_If}=\Ft{1_{I\cap(0,\infty)}f},\quad f\in \cS,\;I\in \cI.$$
$Q_0$ is an idempotent whose image is $G_0(0,\infty)$ and we regard it as 
an operator from $G_0$ to $G_0(0,\infty)$. 
When the restriction of $Q_0$ to $G_0(I)$ for $I\supset (0,\infty)$ is 
bounded, we denote by $Q^I$ the unique bounded extension of $Q_0|_{G_0(I)}$ 
in $\B(G(I),G(0,\infty))$. 
We set $Q=Q^\R$. 

\begin{lemma} Let $M\in \HD_b$ and $T_t=e^{tA_M}$. 
Let $f\in \Hi$ such that $\Lt f(z)=\Lt g(z)M(z)$ with 
$g\in \cS_{(0,\infty)}$. 
Then for every $w\in \rH$ in the resolvent set of $A_M$, we have
$$(T_te_w,f)=\inpr{Fe_w}{Q_0 V_{-t}F\overline{f}}.$$
\end{lemma}

\begin{proof} 
Easy computation yields 
$$(zI-A)^{-1}e_w=\frac{e_w-e_z}{z-w},$$
and so
$$((zI-A)^{-1}e_w,f)=\frac{\Lt g(w)M(w)-\Lt g(z)M(z)}{z-w}.$$
Lemma 3.3 implies  
\begin{eqnarray*}((zI-A_M)^{-1}e_w,f)
&=&(\xi_{M,z},e_w)(e_z,f)+\frac{\Lt g(w)M(w)-\Lt g(z)M(z)}{z-w}\\
&=&M(w)\frac{\Lt g(w)-\Lt g(z)}{z-w}.
\end{eqnarray*}
Since 
\begin{eqnarray*}
\int_0^\infty (T_te_w,f)e^{-tz}dt&=&((zI-A_M)^{-1}e_w,f)=M(w)\frac{\Lt g(w)-\Lt g(z)}{z-w}\\
&=&M(w)\int_0^\infty dt\int_0^\infty g(s+t)e^{-sw-tz},
\end{eqnarray*}
we get $(T_te_w,f)=M(w)\Lt {S_t^*g}(w).$

On the other hand, using 
$Q_0V_{-t}F\overline{f}=Q_0V_{-t}\cF[\overline{g}]=\cF[S_t^*\overline{g}],$
we get 
\begin{eqnarray*}
\inpr{Fe_w}{Q_0 V_{-t}F\overline{f}}_G
&=&\frac{1}{2\pi}\int_{-\infty}^\infty 
\frac{1}{(w+i\lambda)M(i\lambda)}
\overline{\cF[S_t^*\overline{g}](\lambda)}|M(i\lambda)|^2d\lambda\\
&=&\frac{1}{2\pi}\int_{-\infty}^\infty 
\frac{M(-i\lambda)\cF[S_t^*g](-\lambda)}{w+i\lambda}d\lambda\\
&=&\frac{-1}{2\pi}\int_{-\infty}^\infty 
\frac{M(i\lambda)\Lt{S_t^*g}(i\lambda)}{i\lambda-w}d\lambda.
\end{eqnarray*}
Note that $M(z)/(1+z)$ and $\Lt{S_t^*}(z)$ belong to $\Ha$ and 
$$\frac{-1}{2\pi}\int_{-\infty}^\infty 
\frac{M(i\lambda)\Lt{S_t^*g}(i\lambda)}{i\lambda-w}d\lambda
=\lim_{\varepsilon\to+0}\frac{-1}{2\pi i}
\int_{\varepsilon-i\infty}^{\varepsilon+i\infty} 
\frac{M(z)\Lt{S_t^*g}(z)}{z-w}dz.$$
The residue theorem implies that for 
$0<\varepsilon< \RE w<r$ the following holds: 
\begin{eqnarray*}\lefteqn{\frac{-1}{2\pi i}
\int_{\varepsilon-i\infty}^{\varepsilon+i\infty} 
\frac{M(z)\Lt{S_t^*g}(z)}{z-w}dz}\\
&=&M(w)\Lt{S_t^*g}(w)+
\frac{-1}{2\pi i}\int_{r-i\infty}^{r+i\infty}\frac{M(z)\Lt{S_t^*g}(z)}{w-z}dz.
\end{eqnarray*}
Since the second term of the right-hand side tends to 0 as $r$ goes to 
$+\infty$, we get the statement. 
\end{proof}

Let $P$ be the (orthogonal) projection from $G(-\infty,\infty)$ onto 
$G(0,\infty)$. 
We denote by $P^I$ the restriction of $P$ to $G(I)$ for 
$I\supset (0,\infty)$. 

\begin{theorem} Let $M\in \HD$ be an outer function with 
$\overline{M(\overline{z})}=M(z)$. 
\begin{itemize}
\item [$(1)$] The following conditions are equivalent: 
\begin{itemize}
\item [$(\rm{i})$] $M\in \HD_b$. 
\item [$(\rm{ii})$] For all $t>0$, the restriction of $Q_0$ to 
$G_0(-t,\infty)$ is bounded.  
\end{itemize}
When these conditions hold, 
$$Fe^{tA_M}F^{-1}=V_t{Q^{(-t,\infty)}}^*|_{G(0,\infty)},$$
and in particular
$$||e^{tA_M}-S_t||=||Q^{(-t,\infty)}-P^{(-t,\infty)}||.$$
\item [$(2)$]  The following conditions are equivalent:
\begin{itemize}
\item [$(\rm{i})$] $M\in \HD_2$. 
\item [$(\rm{ii})$] For all $t>0$, the restriction of $Q_0$ to 
$G_0(-t,\infty)$ is bounded and $Q^{(-t,\infty)}-P^{(-t,\infty)}$ 
is a Hilbert-Schmidt operator. 
\end{itemize}
When these conditions hold, 
$$||e^{tA_M}-S_t||_{H.S.}=||Q^{(-t,\infty)}-P^{(-t,\infty)}||_{H.S.}.$$
\end{itemize}
\end{theorem}

\begin{proof} (1) Assume $M\in \HD_b$. 
The previous lemma shows that $Q_0$ restricted to $G_0(-t,\infty)$ is 
bounded for all $t>0$ and 
$F{e^{tA_M}}^*F^{-1}=Q^{(-t,\infty)}V_{-t}|_{G(0,\infty)}$. 
Since 
$$FS_tF^{-1}=V_t|_{G(0,\infty)}=PV_tP|_{G(0,\infty)},$$ 
we get $FS_t^*F^{-1}=P^{(-t,\infty)}V_{-t}|_{G(0,\infty)}$ and so 
$||e^{tA_M}-S_t||=||Q^{(-t,\infty)}-P^{(-t,\infty)}||$ holds.

Assume conversely that the restriction of $Q_0$ to $G_0(-t,\infty)$ 
is bounded for all $t>0$ and we set $T_t=F^{-1}V_t{Q^{(-t,\infty)}}^*F$. 
Since $Q^{(-t,\infty)}V_{-t}V_t|_{G_0(0,\infty)}=I$ and 
$$Q^{(-t,\infty)}V_{-t}Q^{(-s,\infty)}V_{-s}|_{G_0(0,\infty)}
=Q^{(-s-t,\infty)}V_{-(s+t)}|_{G_0(0,\infty)},\quad s,t>0,$$
$\{T_t\}_{t\geq 0}$ is a $C_0$-semigroup satisfying the condition (C1). 
The proof of the previous lemma implies 
$$\int_0^\infty(T_te_w,f)e^{-zt}dt=((zI-A_M)^{-1}e_w,f)$$
and so $M\in \HD_b$ and $T_t=e^{tA_M}$. 

(2) follows from (1). 
\end{proof}

The condition in (2) that $Q^{(-t,\infty)}-P^{(-t,\infty)}$ is a 
Hilbert-Schmidt operator is equivalent to that the orthogonal 
projection from $G(0,\infty)$ to $G(-t,0)$ is a Hilbert-Schmidt 
operator. 
Indeed, the following statement holds in general: 

\begin{lemma} Let $H$ be a Hilbert space which is a topological direct sum 
of two closed subspaces $H_1$ and $H_2$ and let $P_i$ be the orthogonal 
projection onto $H_i$ for $i=1,2$. 
Let $Q$ be the (not necessarily orthogonal) projection from $H$ onto $H_1$ 
with respect to the decomposition $H=H_1\oplus H_2$, 
Then the following two conditions are equivalent: 
\begin{itemize}
\item [(1)] The orthogonal projection from $H_1$ to $H_2$ is a 
Hilbert-Schmidt operator. 
\item [(2)] $Q-P_1$ is a Hilbert-Schmidt operator. 
\end{itemize}
\end{lemma}

\begin{proof} \cite[p.308]{Ta} implies that $H$ has an orthogonal 
decomposition $H=K_1\oplus K_2\oplus (K_3\oplus K_3)$ such that 
$$P_1=1\oplus 0\oplus
\left[\begin{array}{cc}
1&0\\
0&0\end{array}\right],$$
$$P_2=0\oplus 1\oplus
\left[\begin{array}{cc}
c^2&cs\\
cs&s^2\end{array}\right],$$
where $c\in \B(K_3)$ is a positive contraction and 
$s\in \B(K_3)$ is a positive invertible contraction with 
$c^2+s^2=I$. 
Then $Q$ is given by 
$$Q=1\oplus 0\oplus
\left[\begin{array}{cc}
1&-s^{-1}c\\
0&0\end{array}\right].$$
Thus we get $||P_2P_1||_{H.S.}^2={\rm Tr}(c^2)$ and 
$||Q-P_1||_{H.S.}^2={\rm Tr}(s^{-2}c^2)$.  
\end{proof}

Note that the following three conditions are equivalent: 
(1) $Q_0$ is bounded. (2) The restriction of 
$Q_0$ to a $Q_0$-invariant dense subspace $G_1$ of $G_0$ is bounded. 
(3) $G=G(-\infty,0)\oplus G(0,\infty)$ topologically. 
As $G_1$, we follow Tsirelson \cite{T2} and adopt the linear span of
$$\{\frac{1}{(1+i\lambda)^m},\frac{1}{(1-i\lambda)^n}\}_{m,n\in \N}.$$
In \cite{T2}, B. Tsirelson observed that questions around $Q_0$ 
are reduced to those of ``past-and-future geometry" via the conformal 
transformation 
$$e^{i\theta}=\frac{i\lambda-1}{i\lambda+1}.$$
Let $N(e^{i\theta})=M(i\lambda)$ and let $\mu(d\theta)=
|N(e^{i\theta})|^2d\theta$, which is a {\it finite} measure on 
the unit circle $\T$. 
Let $\tilde{G}=L^2(\T,\mu)$ and we set 
$\tilde{G}_{+}$ and $\tilde{G}_-$ the closed linear spans of 
$\{e^{in\theta}\}_{n\geq 0}$ and $\{e^{in\theta}\}_{n<0}$ 
respectively. 
As we have 
$$\mu(d\theta)=\frac{2|M(i\lambda)|^2}{1+\lambda^2}d\lambda=
\frac{4\pi}{1+\lambda^2}\nu(d\lambda),$$
we can introduce a unitary $W$ from $\tilde{G}$ onto $G$ by setting  
$$Wf(\lambda)=\frac{2\sqrt{\pi}f(\frac{1-i\lambda}{1+i\lambda})}
{1+i\lambda},$$
which maps $\tilde{G}_{+}$ and $\tilde{G}_-$ onto $G(0,\infty)$ and 
$G(-\infty,0)$ respectively. 
Those measures $\mu$ for which $\tilde{G}$ is a topological direct sum 
of $\tilde{G}_{+}$ and $\tilde{G}_-$ are known as Helson-Szeg\"o  measures 
\cite[Chapter VII,D]{Koo}. 

We get the following two corollaries from Tsirelson's observation. 

\begin{cor} Let $M\in \HD$ be an outer function with 
$\overline{M(\overline{z})}=M(z)$. 
Then the following conditions are equivalent: 
\begin{itemize}
\item [(1)] $M\in \HD_b$ and $\sup_{t>0}\{||e^{tA_M}-S_t||\}<\infty$. 
\item [(2)] The Riesz projection $Q_0$ is bounded. 
\item [(3)] The measure $\mu$ satisfies the Helson-Szeg\"o  condition: 
there are real functions $u,v\in L^\infty(\T)$ with 
$||v||_\infty<\pi/2$ such that $|N(e^{i\theta})|^2=e^{u(\theta)+
\tilde{v}(\theta)}$, 
where $\tilde{v}$ is the conjugate function of $v$. 
\item [(4)] The measure $\mu$ satisfies the Hunt-Muckenhoupt-Wheeden 
condition: there exists a positive constant $c$ such that 
for all interval $I$ in $\T$ the following holds
$$\Big[\int_I|N(e^{i\theta})|^2d\theta\Big]
\Big[\int_I\frac{d\theta}{|N(e^{i\theta})|^2} \Big]\leq c|I|^2.$$
\end{itemize}
\end{cor}

\begin{example}
When $M(z)=(1+z)^\alpha$ with $-1/2<\alpha<1/2$, we have 
$$|N(e^{i\theta})|^2=2^{1-2\alpha}|\sin\frac{\theta}{2}|^{-2\alpha},$$
and so it is easy to show that $M(z)$ satisfies the above (4). 
In Corollary 7.6, we directly show (1) as well. 
\end{example}

\begin{cor} Let $M\in \HD$ with $\overline{M(\overline{z})}=M(z)$. 
Then the following conditions are equivalent: 
\begin{itemize}
\item [$(1)$] $M$ is an outer function, $Q_0$ is bounded and 
$Q-P$ is a Hilbert-Schmidt operator.  
\item [$(2)$] $M$ is an outer function satisfying 
$$\int_{-\infty}^\infty\int_{-\infty}^\infty
\frac{|\log|M(is)|-\log|M(it)||^2}{|s-t|^2}dsdt<\infty.$$
\item [$(3)$]
$$\sup_{x>0}\int_{-\infty}^\infty\int_{-\infty}^\infty
\frac{|M(x+iy)-M(i s)|^2}{[x^2+(y-s)^2]|M(x+iy)|^2}dsdy<\infty.$$
\item [$(4)$]
$$\sup_{x>0}\frac{1}{x}\int_{-\infty}^\infty
\big(\frac{\Pin{|M(i\cdot)|^2}(x+iy)}{|M(x+iy)|^2}-1\big)dy<\infty.$$ 
\item [$(5)$] $M\in \HD_2$ and 
$\sup_{t>0}\{||e^{tA_M}-S_t||_{H.S.}\}<\infty$.  
\end{itemize}
\end{cor}

\begin{proof} 
Tsirelson \cite[Theorem 3.2]{T2} (with Lemma 6.5) shows that (2) implies (1). 
The equivalence of (3),(4), and (5) has already been shown. 
Theorem 6.4 shows that (1) implies (5). 

Assume (1). Since $\mu$ is a finite measure, 
Ibragimov-Solev theorem (see \cite[Proposition 1.8]{T2}) implies that 
$\log |N(e^{i\theta})|^2$ belongs to the Sobolev space $W_2^{1/2}(\T)$ 
and equivalently, the function $M$ satisfies the condition of (2). 

Assume (3), (4) and (5) now. 
(4) implies that $M(z)$ has no zeros on $\rH$ and the 
Blaschke component of $M(z)$ is trivial. 
(3) with the argument of the proof of Theorem 4.5 shows 
that $M$ has a trivial singular inner component as well and 
so $M(z)$ is outer. 
Corollary 6.6 implies that the Riesz projection is bounded and 
$Q$ is well-defined. 
Let $P_{-(t,\infty)}$ be the orthogonal projection from $G$ 
onto $G(-t,\infty)$. 
Note that $P_{-(t,\infty)}$ converges to $I$ in the strong operator 
topology as $t$ tends to $+\infty$.  
(5) implies that there is a positive constant $c$ such that 
$||(Q-P)P_{(-t,\infty)}||_{H.S.}\leq c$ for all $t>0$. 
Since the trace of $\B(G)$ is lower semi-continuous in the 
weak operator topology, we get
$$||Q-P||_{H.S.}^2={\rm Tr}((Q-P)(Q-P)^*)
=\lim_{t\to+\infty}{\rm Tr}((Q-P)P_{(-t,\infty)}(Q-P)^*)\leq c.$$
Thus $Q-P$ is a Hilbert-Schmidt operator and (1) holds. 
\end{proof}

We have seen that $M(z)=\log^\alpha(a+z)$ belongs to $\HD_2$ 
for $a\geq 1$ and $\alpha>0$ in Section 4. 
Using Tsirelson's criterion \cite[Proposition 3.6]{T2}, 
we can actually show that $M$ satisfies the conditions of Corollary 6.7 
if and only if $a>1$ ($\alpha$ can be an arbitrary real number). 
$\log^\alpha(1+z)$ with $\alpha>0$ is a typical example of $M \in \HD_2$ 
without satisfying the conditions of Corollary 6.7. 
This is caused by a zero of $M$ on the imaginary axis. 
Indeed, for every $M\in \HD_2$, the function $M(z)z/(1+z)$ belongs 
to $\HD_2$ and does not satisfy the conditions of Corollary 6.7 
(See Proposition 7.11.)

We end this section with showing that the Tsirelson's 
$E_0$-semigroups constructed from off white noises in \cite{T1} 
actually come from our construction in Theorem 2.3. 
Although it is possible to show the statement by a purely measure theoretical 
argument as in \cite{T1} (in fact, the author first obtained the statement 
in that way), we take an operator theoretical approach using 
Shale's result \cite{Sh} inspired by Bhat and Srinivasan's paper \cite{BS}. 

Let $\rho(x)$ be a real valued measurable function on $\R$ such that 
$$\int_{-\infty}^{+\infty}\frac{e^{\rho(\lambda)}}{1+\lambda^2}<\infty,$$
$$\int_{-\infty}^\infty\int_{-\infty}^\infty
\frac{|\rho(s)-\rho(t)|^2}{|s-t|^2}dsdt<\infty.$$
We set 
$$M(z)=\exp\big[\frac{1}{2\pi}\int_{-\infty}^\infty 
\frac{\lambda z+i}{\lambda+iz}\cdot\frac{\rho(\lambda)
d\lambda}{1+\lambda^2}\big],$$
and set $\{T_t\}_{t\geq 0}$ to be the $C_0$-semigroup corresponding to $M$. 
In the above setting, we have 
$$\nu(d\lambda)=\frac{e^{\rho(\lambda)}}{2\pi}d\lambda.$$
Thus our $\nu$ is nothing but a multiple of $\nu$ in \cite[Section 9]{T1} 
and $G(s,t)$ is the complexification of $G_{s,t}$ in \cite[Section 9]{T1}. 
We set $G(s,t)_\R$ to be the real part of $G(s,t)$, that is, 
$G(s,t)_\R$ is the closure of the linear span of the functions of the form 
$\Ft{1_If}$ where $I\subset (s,t)$ and $f\in \cS$ is real valued. 

Let 
$$L: G(-\infty,0)\oplus G(0,\infty)\ni f\oplus g\mapsto f+g\in G,$$
where we regard $G(-\infty,0)\oplus G(0,\infty)$ as an 
orthogonal direct sum. 
Since $I-L^*L$ is a Hilbert-Schmidt class operator, Shale's result \cite{Sh} 
shows that there exists a unitary operator 
$\Gamma(L):e^{G(-\infty,0)}\otimes e^{G(0,\infty)}\rightarrow e^G$ satisfying 
\begin{eqnarray*}
\Gamma(L)\big(W(f_1+ig_1)\otimes W(f_2+ig_1)\big)\Gamma(L)^*
&=&W(L(f_1\oplus f_2)+i{L^{*}}^{-1}(g_1\oplus g_2))\\
&=&W(f_1+f_2+i{L^{*}}^{-1}(g_1\oplus g_2)),
\end{eqnarray*}
where $f_1,g_1\in G(-\infty,0)_\R$ and $f_2,g_2\in G(0,\infty)_\R$. 
Then Tsirelson's $E_0$-semigroup $\beta_t$ acting on $\B(e^{G(0,\infty)})$ 
is given by 
$$1\otimes \beta_t(W(f+ig))=\mathrm{Ad}(\Gamma(L)\Gamma(V_t)\Gamma(L)^*)
(W(f+ig)).$$

\begin{theorem} Let the notation be as above. 
Then Tsirelson's $E_0$-semigroup $\{\beta_t\}_{t\geq 0}$ is conjugate to 
the $E_0$-semigroup $\{\alpha_t\}_{t\geq 0}$ acting on $\B(e^{\Hi})$ given by 
$$\alpha_t(W(f+ig))=W(S_tf+iT_tg),\quad f,g\in \Hi_\R.$$
\end{theorem}

\begin{proof} Direct computation shows that for $f,g\in G(0,\infty)_\R$, 
we have
\begin{eqnarray*}
1\otimes \beta_t(W(f+ig))&=&
W(L^{-1}V_tf+iL^*V_t{L^{*}}^{-1}(0\oplus g))\\
&=&W((0\oplus V_tf)+iL^*V_t{L^{*}}^{-1}(0\oplus g)).
\end{eqnarray*}
Let $f_1\in G(-\infty,0)_\R$ and $f_2\in G(0,\infty)_\R$. 
Then 
\begin{eqnarray*}
\inpr{L^*V_t{L^{*}}^{-1}(0\oplus g)}{f_1\oplus f_2}
&=&\inpr{0\oplus g}{L^{-1}V_{-t}L(f_1\oplus f_2)}\\
&=&\inpr{0\oplus g}{L^{-1}V_{-t}(f_1+f_2)}\\
&=&\inpr{0\oplus g}{V_{-t}f_1+(I-Q^{(-t,\infty)})V_{-t}f_2\oplus 
Q^{(-t,\infty)}V_{-t}f_2}\\
&=&\inpr{g}{Q^{(-t,\infty)}V_{-t}f_2}.
\end{eqnarray*}
Thus Lemma 6.3 shows the statement. 
\end{proof}
 
\begin{remark} Let $\varphi \in L^1_{\mathrm{loc}}[0,\infty)\cap 
L^2((0,\infty),(1\wedge x)dx)$ with $\|\varphi\|_1<1$ and 
we set $M(z)=1-\Lt{\varphi}(z)$. 
Then there exist positive constants $0<c_1\leq 1\leq c_2$ such that 
$c_1\leq |M(i\lambda)|\leq c_2$ holds for all $\lambda$. 
As a system of topological vector spaces, $\{G(s,t)\}_{s<t}$ and 
$\{L^2(s,t)\}_{s<t}$ are canonically isomorphic. 
Since Tsirelson's infinitesimal sequence invariant in \cite{T1} is 
rather an invariant for $\{G(s,t)\}_{s<t}$ as a system of topological 
vector spaces, it does not distinguish the $E_0$-semigroup arising from 
$M$ from the CCR flow. 
Yet, we will show that even such $M$ sometimes produces 
a $E_0$-semigroup of type III in the forthcoming paper \cite{IS}. 
Namely we will show that the resulting $E_0$-semigroup is of type III 
if and only if $\varphi\notin \Hi$ and that there are uncountably many 
mutually non cocycle conjugate $E_0$-semigroups of type III arising 
in this way. 
\end{remark}

\section{Appendix.}
In this appendix, we show how to compute $K_t=e^{tA_M}-S_t$ for a 
function $M$ in $\HD_b$ that is not necessarily in $\HD_2$ and present 
computation of concrete examples. 

\begin{lemma}
Let $M\in \HD$ and $f\in D(A)$. 
Then 
$$(\xi_{M,z},f)=\frac{1}{M(z)} \int_0^\infty e^{-tz}(q,S_t(f'-f))dt,$$
where $q\in \Hi\setminus D(A)$ with $M(z)=(1+z)\Lt q(z)$. 
\end{lemma}

\begin{proof} Thanks to Equation (3.2), we get
\begin{eqnarray*}M(z)(\xi_{M,z},f)&=&
(q,f)-(1+z)\int_0^\infty e^{-tz}(q,S_tf)dt\\
&=&(q,f)-(1+z)(q,(zI-A)^{-1}f)
=-(q,(zI-A)^{-1}(A+I)f)\\
&=&\int_0^\infty e^{-tz}(q,S_t(f'-f))dt. 
\end{eqnarray*}
\end{proof}

\begin{lemma} Let $M\in \HD$ such that there exist a natural number 
$n$ and positive constants $a$ and $C$ such that 
$|M(z)|\geq C(1+|z|)^{-n}$ for $\RE z>a$.  
\begin{itemize}
\item [$(1)$]
There exists a distribution $k\in \cD'(\R\times (0,\infty))$ with support in 
$[0,\infty)\times (0,\infty)$ 
such that for every $f\in \cD(0,\infty)$, $g\in \cD(\R)$ and $b>a$, 
$$\frac{1}{2\pi i}\int_{b-i\infty}^{b+i\infty}
\big(\int_{-\infty}^\infty e^{xz}g(x)dx(\xi_{M,z},f)\big)dz=
(k,g\otimes f).$$
\item [$(2)$] Moreover, if there exists $r\in \Hi$ such that 
$\Lt r(z)=1/[(1+z)M(z)]$, then $(k,g\otimes f)=(k_f,g)$ where
$$k_f(x)=-(q,f-f')r(x)-\int_0^xr(s)(q,S_{x-s}(f-2f'+f''))ds.$$
\end{itemize}
\end{lemma}

\begin{proof} (1) Let $l(t)=(q,S_t(f'-f))$, which is a bounded smooth 
function on $[0,\infty)$ vanishing at infinity. 
Since there exists positive constants $a$ and $C'>0$ 
such that $$|(\xi_{M,z},f)|\leq C'(1+|z|)^n,\quad \RE z>a,$$
there exists a distribution $k_f\in \cD'(\R)$ with support in $[0,\infty)$ 
such that the integral in the statement (1) is $(k_f,g)$. 
Note that it does not depend on $b>a$. 
By definition, for fixed $f$ the map $\cD(\R)\ni g\mapsto (k_f,g)$ 
is continuous. 
We claim that for fixed $g$, the map $\cD(0,\infty)\ni f\mapsto (k_f,g)$ 
is continuous. 
Note that the function 
$$y\mapsto \int_{-\infty}^\infty e^{x(b+iy)}g(x)dx,$$
is rapidly decreasing.  
Since 
$$|(\xi_{M,b+iy},f)|\leq C'(1+|b+iy|)^n||q||(||f||+||f'||),$$
there exists a constant $C_g$, depending on $g$, 
such that $|(k_f,g)|\leq C_g(||f||+||f'||)$ holds for all 
$f\in \cD(0,\infty)$, which shows the claim. 
Thus the Schwartz nuclear theorem \cite{S} implies the statement. 

(2) Since 
$$(\xi_{M,z},f)=(1+z)\Lt r(z)\Lt l(z)=\Lt r(z)\big(\Lt{l+l'}(z)+l(0)\big),$$
we get $k_f= r*(l+l')+(q,f'-f)r$. 
\end{proof}

\begin{cor} Let $M$, $k$, $q$, and $r$ be as above. 
\begin{itemize}
\item [$(1)$] If the derivatives $q'$ and $q''$ in $\cD'(0,\infty)$ 
belong to $L^2_{\mathrm{loc}}(0,\infty)$, then 
$k$ is in $L^1_{\mathrm{loc}}(0,\infty)^2$ and is given by
\begin{eqnarray*}\lefteqn{k(x,y)=-r(x)\big(q(y)+q'(y)\big)}\\
&-&\int_0^x r(s)\big(q(x+y-s)+2q'(x+y-s)+q''(x+y-s)\big)ds. 
\end{eqnarray*}
\item [$(2)$] If $r'$ and $q'$ in $\cD'(0,\infty)$ belong to 
$L^2_{\mathrm{loc}}(0,\infty)$ and $r$ is continuous at 0, then 
$k$ is in $L^1_{\mathrm{loc}}(0,\infty)^2$ and is given by
\begin{eqnarray*}\lefteqn{k(x,y)=-r(0)\big(q(x+y)+q'(x+y)\big)}\\
&-&\int_0^x \big(r(s)+r'(s)\big)\big(q(x+y-s)+q'(x+y-s)\big)ds. 
\end{eqnarray*}
\end{itemize}
\end{cor}

\begin{theorem} Let $M\in \HD$ with positive constants $n$, $a$, and $C$ 
as in Lemma 7.2 and let $k$ be the distribution defined in 
Lemma 7.2. 
Then $M\in \HD_b$ if and only if there exist positive constants $D$ and $b$ 
such that for all $t\geq 0$, 
$$|(k,g\otimes f)|\leq De^{bt} ||f||_2\cdot||1_{(0,\infty)}g||_2,\quad 
f\in \cD(0,\infty),\; g\in \cD((-\infty,t)).$$
When this assumption is satisfied, the semigroup 
$\{e^{tA_M}=S_t+K_t\}_{t\geq 0}$ is given as follows:  
Let $\tilde{K}_t\in \B(\Hi,L^2(0,t))$ determined by 
$(k,g\otimes f)=(\tilde{K}_tf,1_{(0,t)}g)$. 
Then 
$$K_tf(x)=\left\{
\begin{array}{ll}\tilde{K}_tf(t-x)& (x<t)\\
0&(t\leq x)
\end{array}
\right.$$
\end{theorem}

\begin{proof} When $M\in \HD_b$, necessity of the constants $D$ and $b$ 
with the above property is obvious. 

Assume conversely that there exist $D>0$ and $b$ satisfying 
the conditions above and define $\tilde{K}_t$ and $K_t$ as above. 
The family $\{K_t\}_{t\geq 0}$ is strongly continuous with 
$||K_t||\leq De^{bt}$ and the Laplace transformation of 
$\{K_t\}_{t\geq 0}$ exists for $\RE z>b$. 
Let $f\in \cD(0,\infty)$ and $g\in \cD(\R)$. 
Note that for $0<s<t$, $\tilde{K}_sf(x)=\tilde{K}_tf(x)$ for almost all $x\in (0,t)$.  
We define $Kf(x)$ to be $\tilde{K}_tf(x)$ with $x<t$.   
Note that we have 
$$\int_0^t|Kf(x)|^2dx\leq D^2e^{2bt}||f||^2.$$
Then 
\begin{eqnarray*}\int_0^\infty e^{-tz}(K_tf,g)dt
&=&\int_0^\infty dt e^{-tz}\int_0^t Kf(t-x)g(x)dx\\
&=&\int_0^\infty dx g(x)\int_0^\infty e^{-(t+x)z}Kf(t)dt\\
&=&(g,e_z)\int_0^\infty e^{-tz}Kf(t)dt\\
&=&(g,e_z)\lim_{m\to\infty}(k,h_me_z\otimes f),
\end{eqnarray*}
where $\{h_m\}_{m=0}^\infty$ is an increasing sequence in 
$\cD(-\infty,\infty)$ converging pointwisely to a smooth function $h$ 
satisfying the following properties: 
$0\leq h_m(x)\leq 1$ for all $x\in \R$ and $h_m(x)=1$ for $[-1,m]$ and 
$h_m(x)=0$ for $x\leq -2$. 
We show 
$$\lim_{m\to\infty}(k,h_me_z\otimes f)=(\xi_{M,z},f).$$
Indeed, let $k_f$ be the distribution whose Laplace transformation 
is $(\xi_{M,z},f)$. 
Note that $e_ak_f$ is a tempered distribution. 
Thus for $\RE z>a,b$, we have 
\begin{eqnarray*}\lim_{m\to\infty}(k,h_me_z\otimes f)
&=&\lim_{m\to\infty}(k_f,e_zh_m)=\lim_{m\to\infty}(e_ak_f,e_{z-a}h_m)
=(e_ak_f,e_{z-a}h)\\
&=&\Lt{e_ak_f}(z-a)=\Lt{k_f}(z)=(\xi_{M,z},f). 
\end{eqnarray*}

Let $T_t=S_t+K_t$, whose Laplace transformation is the resolvent of 
$A_M$. 
Since there exists a constant $D_1>0$ such that 
$||T_t||\leq D_1e^{bt}$, we have 
$$||(zI-A_M)^{-m}||\leq \frac{D_1}{(\RE z-b)^m},\quad 
m\in \N,\; \RE z>b,$$
which shows that $A_M$ generates a $C_0$-semigroup thanks to \cite{Y}, 
which should coincides with $\{T_t\}_{t\geq 0}$. 
\end{proof}

When $M(z)=(1+z)^\alpha$ with $-1/2<\alpha<1/2$, we have 
$q(x)=x^{-\alpha}e^{-x}/\Gamma(1-\alpha)$ and 
$r(x)=x^{\alpha}e^{-x}/\Gamma(1+\alpha)$ and so 
$$q(x)+q'(x)=\frac{-\alpha x^{-(1+\alpha)}e^{-x}}{\Gamma(1-\alpha)},$$
$$q(x)+2q'(x)+q''(x)=\frac{\alpha(1+\alpha) x^{-(2+\alpha)}e^{-x}}
{\Gamma(1-\alpha)}.$$
Thus 
\begin{eqnarray*}\lefteqn{k(x,y)
=\frac{\alpha e^{-(x+y)}}{\Gamma(1+\alpha)\Gamma(1-\alpha)}
\Big[x^\alpha y^{-(\alpha+1)}-(\alpha+1)
\int_0^x s^\alpha(x+y-s)^{-(\alpha+2)}ds\Big]}\\
&=&\frac{x^\alpha e^{-(x+y)}}{\Gamma(\alpha)\Gamma(1-\alpha)}
\Big[y^{-(\alpha+1)}-(\alpha+1)x
\int_0^1 t^\alpha(x+y-xt)^{-(\alpha+2)}dt\Big]\\
&=&\frac{\sin (\pi\alpha) x^\alpha e^{-(x+y)}}{\pi}
\Big[y^{-(\alpha+1)}-(\alpha+1)x
\int_1^\infty [(x+y)u-x]^{-(\alpha+2)}du\Big]\\
&=&\frac{\sin (\pi\alpha) x^\alpha e^{-(x+y)}}{\pi}
\Big[y^{-(\alpha+1)}-\frac{xy^{-(\alpha+1)}}{x+y}\Big]\\
&=&\frac{\sin (\pi\alpha)}{\pi}\cdot 
\frac{x^\alpha y^{-\alpha} e^{-(x+y)}}{x+y}.
\end{eqnarray*}

\begin{lemma} Let 
$$k(x,y)=\frac{\sin (\pi\alpha)}{\pi}\cdot 
\frac{x^\alpha y^{-\alpha} e^{-(x+y)}}{x+y}.$$
For $-1/2<\alpha<1/2$, the operator $K$ given by 
$$Kf(x)=\int_0^\infty k(x,y)f(y)dy$$
is a bounded operator in $\B(\Hi)$ with $||K||\leq |\tan(\alpha\pi)|$. 
\end{lemma}

\begin{proof} We apply the Schur test \cite[Theorem 5.2]{HS} to $k$. 
Indeed, using the fact
$$\int_0^\infty\frac{1}{(1+t)t^\beta}dt=\frac{\pi}{\sin(\beta\pi)}, 
\quad 0<\beta<1,$$
we get 
$$\int_0^\infty|k(x,y)|\frac{1}{\sqrt{x}}dx
\leq \frac{|\sin(\alpha\pi)|}{\pi\sqrt{y}}\int_0^\infty 
\frac{dt}{(1+t)t^{1/2-\alpha}}=\frac{|\tan (\alpha\pi)|}{\sqrt{y}},$$
$$\int_0^\infty|k(x,y)|\frac{1}{\sqrt{y}}dy
\leq \frac{|\sin(\alpha\pi)|}{\pi\sqrt{x}}\int_0^\infty 
\frac{dt}{(1+t)t^{1/2+\alpha}}=\frac{|\tan(\alpha\pi)|}{\sqrt{x}},$$
which shows the statement. 
\end{proof}

\begin{cor} Let $M(z)=(1+z)^\alpha$ with $-1/2<\alpha<1/2$. 
Then $M\in \HD_b$ with 
$$e^{tA_M}f(x)=1_{(0,t)}(x)Kf(t-x)+1_{[t,\infty)}(x)f(x-t),$$
$$Kf(x)=\frac{\sin (\alpha\pi)}{\pi}\int_0^\infty  
\frac{x^\alpha y^{-\alpha} e^{-(x+y)}}{x+y}f(y)dy.$$
\end{cor}

\begin{cor} $M(z)=1-re^{-az}$ with $a>0$ and $r\in \C$. 
Then $M\in \HD_b$. 
\end{cor}

\begin{proof} Easy computation yields 
$$(k,g\otimes f)=\sum_{n=1}^\infty r^n\int_{(n-1)a}^{na}g(x)f(na-x)dx,$$
which satisfies the assumption of Theorem 7.4. 
\end{proof}

\begin{prop} Let $M\in \HD_b$ such that $1/\big(zM(z)\big)$ is 
the Laplace transformation of a function in $L^2_{\mathrm{loc}}(0,\infty)$, 
that is, there exists $a>0$ satisfying 
$$\sup_{x\geq a}\int_{-\infty}^\infty\frac{dy}{|x+iy|^2|M(x+iy)|^2}<\infty.$$ 
Let $\beta\in \C$ and $\gamma\in \rH$ with $\RE \beta\geq 0$. 
Then 
$$M_1(z)=\frac{z-\beta}{z+\gamma}M(z)$$
belongs to $\HD_b$. 
Let $k$ and $k_1$ be the distributions defined in Lemma 7.2 for $M$ and $M_1$ 
respectively. Then 
$$k_1(x,y)=k(x,y)+(\beta+\gamma)r_\beta(x)q_\gamma(y),$$ 
where $r_\beta\in L^2_{\mathrm{loc}}(0,\infty)$ and 
$q_\gamma\in \Hi$ are determined by  
$$\Lt{q_\gamma}(z)=\frac{M(z)}{z+\gamma},\quad 
\Lt{r_\beta}(z)=\frac{1}{(z-\beta)M(z)}.$$
In particular, $e^{tA_{M_1}}-e^{tA_M}$ is a rank one operator. 
\end{prop}

\begin{proof} Thanks to Lemma 3.3, we have 
$$\Lt{\xi_{M_1,z}}(w)
=\frac{1-\frac{M_1(w)}{M_1(z)}}{z-w}
=\Lt{\xi_{M,z}}(w)
+\frac{\beta+\gamma}{(z-\beta)M(z)}\cdot\frac{M(w)}{w+\gamma},$$
and so 
$$\xi_{M_1,z}=\xi_{M,z}+\frac{\beta+\gamma}{(z-\beta)M(z)}q_\gamma.$$
Now the statement follows from Theorem 7.4. 
\end{proof}



\end{document}